\documentclass[12pt]{article}
             \newlength{\oldparindent}
             \usepackage{amsmath}
             \usepackage{latexsym}
             \usepackage{t1enc}
             \usepackage[latin1]{inputenc}

             \newtheorem{theoreme}{Theorem}[section]
                \newtheorem{Defi}[theoreme]{Definition}
             \newtheorem{proposition}[theoreme]{Proposition}
             \newtheorem{corollaire}[theoreme]{Corollary}
             \newtheorem{lemme}[theoreme]{Lemma}
             \newtheorem{remarque}[theoreme]{Remark}
             \newtheorem{Not}[theoreme]{Notation}
             \newtheorem{hyp}[theoreme]{Assumption}

             \font\myf=msbm10 at 12pt
             \newcommand{\R}{\hbox{\myf R}}
             \newcommand{\1}{{\bf 1}}
             \newcommand{\pr}{\hbox{\myf P}}
             \newcommand{\esp}{\hbox{\myf E}}
             
             \newcommand{\N}{\hbox{\myf N}}
             \newcommand{\C}{\hbox{\myf C}}

             \newcommand{\F}{{\cal F}}

             \newcommand{\eps}{\varepsilon}

\begin{document}

             \def\second {\vtop{\baselineskip=11pt
             \hbox to 125truept{\hss Diana DOROBANTU \hss}\vskip2truept
             \hbox to 125truept{\small\hss U.M.R. CNRS C 5583\hss}
             \hbox to 125truept{\small\hss Laboratoire de Statistique et Probabilit\'es\hss}
             \hbox to 125truept{\small\hss Universit\'e Paul Sabatier\hss}
             \hbox to 125truept{\small\hss 118 route de Narbonne\hss}
             \hbox to 125truept{\small\hss  31 062 TOULOUSE cedex 04\hss}
             \hbox to 125truept{\small\hss FRANCE\hss}
             \hbox to 125truept{\small\hss dorobant@cict.fr\hss}}}

            \title{{ Optimal stopping  for Lévy processes and affine functions  }}
             \author{ \second}
             \date{}
             \maketitle
             \vspace{1cm}
             \hbox to 125truept{\hss Abstract \hss}\vskip2truept

 This paper studies an optimal stopping problem for  Lévy processes.   We give a justification of  the form of the Snell envelope using  standard results of optimal stopping. We also justify the convexity  of the value function, and without  a priori restriction to a particular class of stopping times, we deduce that the   smallest  optimal stopping time is necessarily a hitting time. We propose a method  which allows to obtain the optimal threshold.   Moreover this method allows to avoid long calculations of the  integro-differential operator used in the usual proofs.

 \vspace{1cm}
 Keywords : Lévy processes, optimal stopping problem, hitting time, Snell envelope.
\newpage


             \section{Introduction}
             \label {section1}

            In this paper we study an optimal stopping problem for  jump processes and its application in Finance. We choose to solve a particular optimal stopping problem  for Lévy processes. Without a priori restriction  to a particular class of stopping times like in \cite{Chen}, we propose a method to find the optimal stopping time form (it will be a hitting time), as well as for the calculation of the optimal threshold.

            In fact we seek to control a stochastic process V of the form $V=ve^X$ where $v$ is a real strictly positive constant and $X$ a  Lévy process. We consider the following optimal stopping problem :
\begin{equation}
\label{eq1intro}
esssup_{ \tau\in \Delta, \tau\geq t} \esp\left(\int_t^{\tau}e^{-r(s-t)}h(V_s)ds\mid \F^V_t\right),
\end{equation}
 where   $r>0$,  $\F^V_t=\sigma(V_s, ~s\leq t)$,  $\Delta$ is the set of   
$\F^V_.$-stopping times and $h$ is an affine function.   We will be brought back  to find a stopping time  $\tau^*$ which maximizes $\tau\mapsto \esp_v(A_{\tau})$ where $A_{\tau}=e^{-r\tau}f(V_{\tau})$ and $f$ an affine function. In many papers the optimal stopping time is supposed from the beginning to be a hitting time, here we show that the optimal stopping time is necessarily of the form $\tau_b=inf\{t\geq 0 : ~~V_t\leq b\}$. Following \cite{Sh}, we also introduce  a decreasing sequence of almost surely finite stopping times  $(\tau_{\eps}, ~\eps>0)$ which converges to the  optimal stopping time ; this is about  the   $\eps$-optimal stopping times (i.e. $\esp(A_{\tau^*}) -\eps \leq \esp(A_{\tau_{\eps}} ))$.  We give a justification of the form of the Snell envelope of the process $A$ using  standard results of optimal stopping of \cite{KLM, Sh}, we argue the convexity of the function $v\mapsto \esp(A_{\tau^*}|V_0=v)$ and  the optimal stopping time form. The main result is given by Theorems \ref{th1} and \ref{th2} which allow to determine the optimal threshold. The method used here allows to solve the optimal stopping problem when  the joint  Laplace transform of  $(\tau_b, ~X_{\tau_b})$, i.e. $\esp[e^{-r\tau_b+aX_{\tau_b}}]$, is known. 

The optimal stopping theory  is a subject which often appears in the specialized literature, having applications for example in Medicine \cite{medicine} or Finance \cite{KS2}.  Among others, Leland \cite{Leland1, Leland2, Leland3}, Duffie and Lando \cite{Duf} or Villeneuve \cite{Vil} studied the optimal stopping problem for a diffusion process. Moreover, there are other authors who used mixed  diffusion-jumps processes for their models. For example,  Hilberink and Rogers \cite{Hil} or Kyprianou \cite{Kyp} use a spectrally negative Lévy process and Le Courtois and Quittard-Pinon \cite{Courtois} stable Lévy processes. Mixed diffusion-jump processes with double   exponential jumps  were studied by Chen and Kou \cite{Chen}, Kou and Wang \cite{Kou, Kou2}, Dao \cite{Dao}.  In \cite{Mord}, the jumps follow an exponential law and they are either all positive or all negative. In \cite{Dao}, Dao studies a model where the jumps follow an  uniform law. In \cite{pham}, Pham  uses a jump-diffusion process and the jumps are not restricted to any particular law.

This paper is organized as follows : we introduce the optimal stopping problem
(Section \ref{section2_ch1}). The following section (Section \ref{section3.1_ch1}) contains the main results which characterize the optimal stopping time and the optimal threshold. At the end of this paper we solve   the studied optimal stopping problem in the case of some particular Lévy processes : Brownian motion, Poisson process, double exponential jump-diffusion process, a particular Lévy process with  positive  jumps on the one hand (Section \ref{section4_ch1}), spectrally negative Lévy processes on the other hand (Section \ref{levneg}). We will recover Duffie and Lando's result \cite{Duf} for a Brownian motion with drift and Kou and Wang's result  \cite{Kou2} for a mixed diffusion-jump process with double   exponential jumps. Section \ref{appendices} contains some optimal stopping tools.

             \section{Optimal stopping problem}
              \label {section2_ch1}
              
                Let $V$ be a stochastic process on a filtered probability space $(\Omega,\F,(\F_t)_{t\geq 0},\pr)$. Assume that $V$ is of the form  $$V=ve^X$$ where $v$ is  a real strictly positive constant and $X$ is a Lévy process. 
We  sometimes use the notation $V^v=ve^X$, for $v>0$.

             Following Lévy-Khitchine formula (see for example \cite{Bert1} or \cite{Kyp}), the characteristic function of  $X$ is 
             \begin {equation}
             \nonumber
             \esp(e^{i\lambda X_t})=e^{-t\Psi(\lambda)}
             \end {equation}
             where $\lambda\in \R$ and the function $\Psi : \R \rightarrow \C$ has the form
             \begin {equation}
             \nonumber
             \Psi(\lambda)=-im\lambda+\frac{\sigma^2}{2}\lambda^2+\int_{\R}(1-e^{i\lambda 
             x}+i\lambda x\1_{|x|<1})\Pi(dx)
             \end {equation}
             with $m\in\R$, $\sigma>0$  and $\Pi$ a mesure on $\R^*$ such that $\int 
             (1\wedge |x|^2)\Pi(dx)<\infty.$
             
             From now on, $\esp(.|V_0=v)$ and $\pr(.|V_0=v)$ are denoted $\esp_v(.)$ and $\pr_v(.)$.

             \begin{hyp}
             \label{h1}
             $\esp_v(V_t)<\infty$ for $t\geq 0$.
             \end{hyp}
             
            The condition  $\esp_v(V_t)<\infty$ for $t\geq 0$ is equivalent to $\esp(e^{X_t})<\infty$ and, using Theorem   3.6 page 76  of \cite{Kyp}, it is still equivalent to the  condition $\int_{|x|\geq 1}e^x\Pi(dx)<\infty$.  Moreover $\esp(e^{X_t})$ is of the form $\esp(e^{X_t})=e^{t\psi(1)}$ and $\esp_v(V_t)=ve^{t\psi(1)}.$

             Let $\F^V$ be  the right-continuous complete filtration generated by the process $V$, \\$\F^V_t=\sigma(V_s, ~s\leq t)$. We consider the following  optimal stopping  problem :

             \begin {equation}
             \label{2}
             S_t=esssup_{\tau\in \Delta, \tau\geq t} 
             \esp\left(\int_t^{\tau}e^{-r(s-t)}(\alpha V_s-c)ds\mid \F^V_t\right),
             \end {equation}
             where  $r>0$, $\alpha>0$, $c>0$ and $\Delta$ is the set of  $\F^V_.$-stopping times.
     
       \begin{Defi}
             A stopping time $\tau^*_t$ is said to be optimal at time $t$ if it maximizes (\ref{2}), i.e.
             
             \begin{center}
             $\esp\left[\int_t^{\tau^*_t}e^{-r(s-t)}(\alpha V_s-c)ds\mid 
             \F^V_t\right]=esssup_{\tau\in \Delta, \tau\geq t} 
             \esp\left[\int_t^{\tau}e^{-r(s-t)}(\alpha V_s-c)ds\mid 
             \F^V_t\right].$
             \end{center}
             \end{Defi}

             Remark that for every $t\geq 0$, $S_t\geq 0$ because $\tau=t \in \Delta$.
             
             The same type of problem (\ref{2}) was studied by Duffie and Lando in \cite{Duf}. In \cite{Duf}, $X$ is a  Brownian motion with drift.  The authors solve the problem using the Hamilton-Jacobi-Bellman equations.
        
             \begin{hyp}
             \label{h2}
             $r>\psi(1)$.
             \end{hyp}
             The necessity for this assumption is clearly apparent. If Assumption  \ref{h2} were not checked, then
             $\esp_v\left(\int_0^{\infty}e^{-rs}(\alpha V_s-c)ds\right)=\alpha 
             v\int_0^{\infty}e^{-s(r-\psi(1))}ds-\frac{c}{r}$
              would be infinite and $\tau^*_0=\infty$. It implies that the process
             $s\mapsto e^{-rs}(\alpha V_s-c)$ belongs to $L^1(\Omega\otimes 
             \R_+, ~d\pr\otimes ds)$.


             \section{Optimal stopping time}
            \label{section3.1_ch1}
            
            In this part we show that the problem (\ref{2}) admits at least an optimal stopping time and that the smallest one is a hitting time. The proof of this result requires several lemmas.

             \begin{lemme}
             \label{lem1}
            Under Assumptions \ref{h1} and \ref{h2}, for every $\tau\in \Delta$ the following equality is true :
             $$\esp\left(\int_{\tau}^{\infty}e^{-r(s-\tau)}(\alpha V_s-c)ds\mid 
             \F^V_{\tau}\right)=
             \left(\frac{\alpha 
             V_{\tau}}{r-\psi(1)}-\frac{c}{r}\right)1_{\{\tau<\infty\}}.$$
             \end{lemme}

             \begin{preuve}
             \\
             Let $s\geq 0$. The exponential form of the process $V_.$ allows the factorization
             $$V_s=V_{\tau}e^{X_s-X_{\tau}} \hbox{~on the set~}\{s>\tau\}.$$ 
             
             However $X$ is a Lévy process, therefore  $X_s-X_{\tau}$  is independent of 
             $\F^V_{\tau}$ and equal in distribution with $X_{s-\tau}$ conditionally to  $\{s>\tau\}$. Thus,
             $$\esp\left(\int_{\tau}^{\infty}e^{-r(s-\tau)}(\alpha V_s-c)ds\mid 
             \F^V_{\tau}\right)=
             \1_{\{\tau<\infty\}}e^{r\tau}\left[\alpha V_{\tau} 
             \esp\left(\int_{\tau}^{\infty}e^{-rs}e^{X_s-X_{\tau}}ds \mid 
             \F^V_\tau\right)-\frac{e^{-\tau r}c}{r}\right],$$
  from which the result follows :
             $$\esp\left(\int_{\tau}^{\infty}e^{-r(s-\tau)}(\alpha V_s-c)ds\mid 
             \F^V_{\tau}\right)=\left(\frac{\alpha 
             V_{\tau}}{r-\psi(1)}-\frac{c}{r}\right)1_{\{\tau<\infty\}}.$$
             \end {preuve}
             \hfill  $\Box$
             \\

             Using Lemma ~\ref{lem1}, $S_t$ can be rewritten as
     \begin {equation}
             \label{3}
             S_t=\frac{\alpha V_t}{r-\psi(1)}-\frac{c}{r}+e^{rt}
             esssup_{\tau\in \Delta , \tau \geq t} \esp\left[e^{-r\tau}
             \left(\frac{-\alpha 
             V_{\tau}}{r-\psi(1)}+\frac{c}{r}\right)1_{\tau<\infty}\mid\F^V_t\right]
             \end {equation}
             for each $t\geq 0$.
             
             We have to solve the following optimal stopping problem :
             \begin{equation}
             \label{4}
             J_t=esssup_{\tau\in \Delta , \tau \geq t} \esp\left[e^{-r\tau}
             \left(\frac{-\alpha 
             V_{\tau}}{r-\psi(1)}+\frac{c}{r}\right)1_{\tau<\infty}\mid\F^V_t\right].
             \end{equation}

             We introduce the process $Y_.$ defined by  :
             \begin{Not}
             $Y : t\mapsto Y_t=e^{-rt}
             \left(\frac{-\alpha V_t}{r-\psi(1)}+\frac{c}{r}\right).$
             \end{Not}

             \begin{lemme}
             \label{lem2}
              Under Assumptions \ref{h1} and \ref{h2}, the process $Y_.$ converges in  $L^1$ and almost surely and its limit is $Y_\infty=0$.
             \end{lemme}
             \begin{preuve}
             \\
             Since $\esp_v(\mid Y_t\mid)\leq \frac{\alpha 
             v}{r-\psi(1)}e^{-(r-\psi(1))t}+\frac{ce^{-rt}}{r}$, then $Y_t \longrightarrow 
             ^{L^1} 0$. 
             \\
                This process can be written in the form
             $$Y_t=-e^{-(r-\psi(1))t}M_t+N_t, ~t\geq 0$$ where $M$ defined by $M_t=\frac{e^{-\psi(1)t}\alpha 
             V_t}{r-\psi(1)}, ~t\geq 0$ is a positive martingale and $N$ defined by $N_t=\frac{ce^{-rt}}{r}, ~t\geq 0$ is a continuous decreasing bounded  positive function. Consequently, the process $Y_.$ is the difference between a continuous deterministic function which goes to 0 and a positive supermartingale (thus which  converges almost surely).
 Then, the process $Y_.$ converges almost surely when $t$ goes to $\infty$. Moreover the limit of  $Y$ in $L^1$ is equal to $0$, therefore $Y_\infty=0$ almost surely.
             \end {preuve}
             \hfill  $\Box$
             \\
             
             We thus look for an optimal stopping time among almost surely finite stopping times. The random variable 
 $Y_\infty$ being null almost surely,  we can remove the indicator $1_{\tau<\infty}$  in $(\ref{4})$.
              \begin{remarque}
             \label{cpositif}
            We imposed   $c>0$ to avoid the case $c=0$. Let us notice that if $c=~0$, we have to calculate  essential supremum of a negative quantity : 
             $$esssup_{\tau\in \Delta , \tau \geq t} \esp\left[e^{-r\tau}
             \frac{-\alpha V_{\tau}}{r-\psi(1)}\mid\F^V_t\right].$$ In this case for each  $t\geq 0$, the optimal stopping time is $\tau^*_t=\infty$, $J_t=0$ and the optimal value is $S_t=  \frac{\alpha V_t}{r-\psi(1)}.$
             
             In all the particular cases which we will study,  the optimal stopping time goes to infinity when $c$ goes to $0$.
             \end{remarque}
             
             We suppose that the process $Y_.$ checks the following  assumption :
    \begin{hyp}
             \label{h3_ch1}
             The process $Y_.$ is of class $D$ (i.e. the set of random variables $Y_{\tau}$, $\tau \in  \Delta$ is uniformly integrable). 
          \end{hyp}
             
             Using Theorem \ref{classeD_DM} of Section \ref{appendix1}, we prove easily the following result :
          
             \begin{lemme}
             \label{classeD}
              Under Assumptions \ref{h1} and \ref{h2}, a sufficient condition for Assumption  \ref{h3_ch1} is 
              \begin{equation}
                \label{limRn}
                lim_{n\rightarrow\infty}\esp\left(e^{-rR_n+X_{R_n}}\1_{R_n<\infty}|X_0=0\right)=0
                \end{equation}
                 where $R_n=inf\{t\geq 0 : e^{-rt+X_t}\geq n\}.$
             \end{lemme}

          This is the condition  we will check in all the examples of Sections \ref{section4_ch1} and \ref{levneg}.

The process $(t\mapsto Y_t, ~t\geq 0)$ being of class $D$, we can apply the results of optimal stopping (see Section \ref{appendix1}). According to Theorem ~\ref{formeSnell}, the Snell envelope $J$ of  $Y$ is of the form
 $\left(e^{-rt}s(V_t)\right)_{t\geq 0}$ (with $J_\infty=0$ because   $Y_\infty=0$). We denote $f(v)=\frac{-\alpha v}{r-\psi(1)}+\frac{c}{r}$ ; Definition  (\ref{3}) gives 
 \begin{equation}
 \label{nouvelleforme}
 S_t=-f(V_t)+e^{-rt}J_t=-f(V_t)+s(V_t), ~t\geq 0.
 \end{equation} Thus the process $(S_t)_{t\geq 0}$ is of the form $(w(V_t))_{t\geq 0}$ where $w$ is a positive Borelian function. Since $\sigma>0$,   the support of  $V_t$ is $\R^*_+$  (as a consequence of Theorem 24.10 i) page 152  of \cite{sato}), so for $t=0$ Definition  (\ref{nouvelleforme}) yields  
              $S_0=-f(v)+s(v)$  and the function $w$  coincides with the function 
               \begin{equation}
               \nonumber
               v\mapsto -f(v)+s(v).
               \end{equation}

 The function $s$ is a (decreasing) convex function because it is the sup of (decreasing) affine functions  :
             $$s(v)=sup_{\tau\geq 0}\esp_v\left[e^{-r\tau}
             \left(\frac{-\alpha 
             V_{\tau}}{r-\psi(1)}+\frac{c}{r}\right)\right]=sup_{\tau\geq 0}
             \esp_1\left[e^{-r\tau}
             \left(\frac{-\alpha 
             vV^1_{\tau}}{r-\psi(1)}+\frac{c}{r}\right)\right].$$
             \begin{remarque}
             The function $s$ being convex, it is thus  continuous.
             \end{remarque}

             Remark that $s$ is a positive function because
             \\
             \\
             $s(v)\geq sup_{t\geq 0} \esp_v\left[e^{-rt}
             \left(\frac{-\alpha V_{t}}{r-\psi(1)}+\frac{c}{r}\right)\right]\geq 
             sup_{t\geq 0} \esp_v\left[e^{-rt}
             \frac{-\alpha V_{t}}{r-\psi(1)}\right]=sup_{t\geq 0}
             \frac{-\alpha ve^{-(r-\psi(1))t}}{r-\psi(1)}=~0.$

Since   $Y_t \longrightarrow ^{p.s.} 0$,  we have the following result which is already shown in  \cite{Vil} where the process is not a Lévy, but a diffusion process. In \cite{Vil}, the author uses the process  trajectories continuity  to show his result.  Since our process is c\`adl\`ag, we have to remake the proof.
             \begin{lemme}
             \label{lem3}
             For $v>0$, let
             $$s(v)=sup_{\tau\geq 0}\esp_v\left[e^{-r\tau}
             \left(\frac{-\alpha 
             V_{\tau}}{r-\psi(1)}+\frac{c}{r}\right)\right]~~\hbox{and}~~s^+(v)=sup_{\tau\geq 0}\esp_v\left[e^{-r\tau}
             \left(\frac{-\alpha 
             V_{\tau}}{r-\psi(1)}+\frac{c}{r}\right)^+\right].$$
             If $\sigma>0$, then under Assumptions \ref{h1}, \ref{h2} and \ref{h3_ch1},  $s^+(v)>0$ and $s(v)=s^+(v)$ for every $v>0$.
             \end{lemme}
             \begin{preuve}
             \\
             We show that if there exists $v_0>0$ such that $s(v_0)<s^+(v_0)$, then there exists $v_1>0$ such that $s^+(v_1)=0$. We prove that this last relation  can not be satisfied.
             
             By construction, for each $v>0$, $s(v)\leq s^+(v)$.
             Let us suppose that there exists $v_0>0$ such that $s(v_0)<s^+(v_0)$.
             
      The process $V_.$ is a right continuous one, the process $Y^+ : t \rightarrow Y^+_t=e^{-rt}\left(\frac{-\alpha V_t}{r-\psi(1)}+\frac{c}{r}\right)^+$ takes its values in $[0, ~\frac{c}{r}]$, then the assumptions of Theorem  ~\ref{th5} of Section \ref{appendix1} are checked for $Y^+$.
We define the function $f^+(v)=\left(\frac{-\alpha  v}{r-\psi(1)}+\frac{c}{r}\right)^+$ ;  the stopping time
              $$\tau^+=inf\{u\geq 0 : f^+(V^{v_0}_u)=s^+(V^{v_0}_u)\}$$ is the smallest optimal stopping time of the problem $$s^+(v_0)=sup_{\tau\geq 0}\esp_{v_0}\left[e^{-r\tau}\left(\frac{-\alpha V_{\tau}}{r-\psi(1)}+\frac{c}{r}\right)^+\right].$$
              
              Since $Y$ converges almost surely to $0$, then $Y^+$ converges almost surely to $0$ and :
              \begin{equation}
              \nonumber
             s^+(v_0)=\esp_{v_0}\left[e^{-r\tau^+}
             \left(\frac{-\alpha 
             V_{\tau^+}}{r-\psi(1)}+\frac{c}{r}\right)^+\right]=\esp_{v_0}\left[e^{-r\tau^+}
             \left(\frac{-\alpha 
             V_{\tau^+}}{r-\psi(1)}+\frac{c}{r}\right)^+\1_{\tau^+<\infty}\right].
             \end{equation}
             
             Using the definition of  $s$ and $s^+$ :
             $$\esp_{v_0}\left[ e^{-r\tau^+}f(V_{\tau^+})\right]\leq s(v_0)<s^+(v_0)=\esp_{v_0}\left[ e^{-r\tau^+}f^+(V_{\tau^+})\right]$$
             and consequently $\esp_{v_0}\left[ e^{-r\tau^+}\left(f(V_{\tau^+})-f^+(V_{\tau^+})\right)\right]<0$, 
             $\pr_{v_0}\left(\{\omega : f(V_{\tau^+})<0\}\right)>0$
             and \\$\pr_{v_0}\left(\{\omega : s^+(V_{\tau^+})=0\}\right)>0.$ 
             
            Thus there exists $v_1$ such that $s^+(v_1)=0$. Then for any stopping time $\tau$, $\pr_{v_1}$-almost surely $e^{-r\tau}f^+(V_{\tau})=0$ and  in particular for every  $t\in \R_+$, $f^+(V_t)=0$. This involves that $\pr_{v_1}$-almost surely $V_t\geq \frac{c(r-\psi(1))}{\alpha r}$ which is a contradiction because the support of $V_t$ is $\R^*_+$ when $\sigma>0$. So $s^+(v)>0$ for every $v\in\R^*_+$ and $s(v)=s^+(v)$.
             \end {preuve}
             \hfill  $\Box$
             \\

     \begin{remarque}
             We have the equality $s(v)=s^+(v)>0$ for every $v>0$. 
             \end{remarque}
             
                   \begin{proposition}
             \label{formetau}
                If $\sigma>0$, then under Assumptions \ref{h1}, \ref{h2} and \ref{h3_ch1}, there exists at least an optimal stopping time for the problem (\ref{4}).
                
                 For any $c>0$
             there exists $b_c>0$  such that the smallest optimal stopping time has the following form 
             $$\tau_{b_c}=inf \{t\geq 0 : ~V_t \leq b_c\}.$$
             \end{proposition}

             \begin{preuve} 
             \\
             Using Lemma \ref{lem3}, the problem (\ref{4}) can be written as $sup_{\tau\geq 0}\esp(Y_{\tau }^+)$. The assumptions of Theorem ~\ref{th5} of Section \ref{appendix1} are checked and the stopping time
               $$\tau^*=inf\{u\geq 0 :
              f^+(V_u)=s^+(V_u)\}$$ is the smallest optimal stopping time.
However $s(v)=s^+(v)>0$ for all $v>0$, so 
              $$\tau^*=inf\{u\geq 0 :
              f(V_u)=s(V_u)\}$$ is the smallest optimal stopping time. 
               The function  $s$ is upper bounded by $\frac{c}{r}$ because $Y_.^+$ is upper bounded by
               $\frac{c}{r}$ and  $lim_{v\downarrow 0} s(v)=lim_{v\downarrow 0} 
             f(v)=\frac{c}{r}$.
             
              Since  $s$ is convex and  $f$ affine, then
             $inf \{v>0 : f(v)<s(v)\}$ is equal to\\$sup \{v>0 : f(v)=s(v)\}$,  and we denote it $b_c$. Indeed, let $b_c'=sup \{v : f(v)=s(v)\}$  and $b_c=inf \{v : f(v)<s(v)\}$. Since $lim_{v\downarrow 0} s(v)=lim_{v\downarrow 0}  f(v)$, then $b_c'$ exists and $b_c'\geq 0$. If $b_c=0$, then $b_c'=0.$
             
             If  $b_c>0$, then for every $v<b_c$, $f(v)=s(v)$ ; in particular $f(b_c-\frac{1}{n})=s(b_c-\frac{1}{n})$.  When $n$ goes to infinity, since $s$ and $f$ are continuous, then $f(b_c)=s(b_c)$, so $b_c\leq b_c'$.  Let us suppose that $b_c<b_c'$, thus there exists $v$, $b_c<v<b_c'$ such that $f(v)<s(v)$. However $s$ is convex :
             $$\frac{s(v)-s(b_c)}{v-b_c}\leq \frac{s(b_c')-s(v)}{b_c'-v}.$$
            Since, by continuity, $f(b_c')=s(b_c')$, then 
            $$\frac{s(v)-f(b_c)}{v-b_c}\leq \frac{f(b_c')-s(v)}{b_c'-v}$$
             Since $s(v)>f(v)$, then 
             $$\frac{f(v)-f(b_c)}{v-b_c}<\frac{s(v)-f(b_c)}{v-b_c}\leq \frac{f(b_c')-s(v)}{b_c'-v}<\frac{f(b_c')-f(v)}{b_c'-v},$$
        which is a contradiction because, since $f$ is affine, then  $\frac{f(v)-f(b_c)}{v-b_c}=\frac{f(b_c')-f(v)}{b_c'-v}=\frac{-\alpha}{r-\psi(1)}$. Consequently $b_c=b_c'$.

             This means that the smallest optimal stopping time $\tau^*$ is also the first entrance time  in  
             $]0, ~b_c]$.
             \end {preuve}
             \hfill  $\Box$
             \\

The smallest optimal stopping time of (\ref{4}) depends on $c$ and from now on we use the notation $\tau^*(c)$. As a consequence of Theorem ~\ref{toptim} of Section \ref{appendix1}, we have that
              \begin{equation}
              \nonumber
              (t\mapsto e^{-r(t\wedge \tau^*(c))}s(V_{t\wedge \tau^*(c)}), ~t\geq 
             0)\mbox{
              is a martingale and }Y_{\tau^*(c)}=e^{-r\tau^*(c)}s(V_{\tau^*(c)}).
              \end{equation}
              
               We introduce an auxiliary function :
              
                \begin{Defi}
               \label{fc_g}
                Let $g : \R_+^*\times ]0, \frac{(r-\psi(1))c}{r\alpha}[ \rightarrow \R_+^*$ be the function defined by 
              \begin {equation}
             \nonumber
             g(v,b)=\esp_v\left[e^{-r\tau_{b}}
             \left(\frac{-\alpha 
             V_{\tau_{b}}}{r-\psi(1)}+\frac{c}{r}\right)\right]
        \end{equation}
        where $\tau_b=inf\{t\geq 0 : V_t\leq b\}$.
             \end{Defi}
             
             If $b\in\R_+$, then $g$ is not necessarily positive. The condition $b\in ]0, \frac{(r-\psi(1))c}{r\alpha}[$ implies the positivity of $g$.
             \begin{remarque}
             \label{intervalle_ch1}
          Under  the assumptions of Proposition \ref{formetau},  there exists $B_c$ such that $g(., B_c)=s(.)$.
  \end{remarque}

Remark that we can explicit $g$ as a function of Laplace transforms
\begin{equation}
              \nonumber
              {\cal{L}}(x)=\esp\left[e^{-r\bar{\tau}_x}|X_0=0\right],~~~  {\cal{G}}(x)=\esp\left[e^{-r\bar{\tau}_x+X_{\bar{\tau}_x}}|X_0=0\right]
              \end{equation}
              where $\bar{\tau}_x=inf\{t\geq 0 : X_t\leq x\}$. Indeed, the function $g$ can be written as
              $$g(v,b)=\frac{-\alpha v}{r-\psi(1)}{\cal{G}}\left(ln\frac{b}{v}\right)+\frac{c}{r}{\cal{L}}\left(ln\frac{b}{v}\right).$$
           
              Now, the aim is  to calculate, when that is possible, the value of the optimal threshold $B_c$ as a function of $\alpha$, $c$, $r$, $\psi(1)$ and the functions ${\cal{L}}$ and ${\cal{G}}$.
              
              Remark that  ${\cal{L}}(x)={\cal{G}}(x)=1$ for $x\geq 0$ and $0\leq {\cal{G}}(x)\leq{\cal{L}}(x)\leq 1$. Moreover these functions  are increasing. Indeed for every $x\leq y<0,$ $\bar{\tau}_y\leq \bar{\tau}_x$, thus ${\cal{L}}$ is an increasing function. Moreover
              $${\cal{G}}(y)=\esp\left[e^{-r\bar{\tau}_y+X_{\bar{\tau}_y}}\1_{X_{\bar{\tau}_y}\leq x}|X_0=0\right]+ \esp\left[e^{-r\bar{\tau}_y+X_{\bar{\tau}_y}}\1_{x<X_{\bar{\tau}_y}\leq y}|X_0=0\right].$$
              
              On $\{X_{\bar{\tau}_y}\leq x\}$,  $\bar{\tau}_y=\bar{\tau}_x$ $~\pr(.|X_0=0)$-almost surely and $X_{\bar{\tau}_y}=X_{\bar{\tau}_x}$.
             On $\{x<X_{\bar{\tau}_y}\leq y\}$, \\$-r\bar{\tau}_y+X_{\bar{\tau}_y}\geq -r\bar{\tau}_x+x\geq -r\bar{\tau}_x+X_{\bar{\tau}_x}$ and the result follows ${\cal{G}}(y)\geq {\cal{G}}(x)$.

              When $\cal{G}$ is discontinuous at $x=0$, $B_c$ is easy to obtain.
              
  \begin{theoreme}
             \label{th2}
          Let $\sigma>0$. Under Assumptions \ref{h1}, \ref{h2} and \ref{h3_ch1}, we suppose that the function 
             $\cal{G}$  is discontinuous at $x=0$.  
             Then the smallest optimal stopping time is $$\tau^*(c)=inf\{t\geq 0 : V_t\leq 
             B_c\},$$ where $B_c=\frac{c(r-\psi(1))}{r\alpha}lim_{x\uparrow 0}\frac{1-{\cal{L}}(x)}{1-{\cal{G}}(x)}$.
             \end{theoreme}
             \begin{preuve}
             \\
             Let $b\in]0, \frac{(r-\psi(1))c}{r\alpha}[$. The function $g$ has the form 
             
                 \begin {equation}
             \nonumber
             g(v, b)=\left\{
             \begin {array}{l@{\quad}l}
             -\frac{\alpha v}{r-\psi(1)}+ \frac{c}{r}& if \ v\leq b\\
         \frac{-\alpha v}{r-\psi(1)}{\cal{G}}\left(ln\frac{b}{v}\right)+\frac{c}{r}{\cal{L}}\left(ln\frac{b}{v}\right) & if \ v>b. \\
             \end{array}
             \right.
             \end{equation}

             If the function $g(.,b)$ is continuous at $b$, then $b$ is solution of 
             \begin{equation}
             \label{bsoleq}
              -\frac{\alpha b}{r-\psi(1)}+ \frac{c}{r}= \frac{-\alpha b}{r-\psi(1)}{\cal{G}}(0^-)+\frac{c}{r}{\cal{L}}(0^-).
              \end{equation}
              However, $\cal{G}$ is discontinuous at $x=0$, so ${\cal{G}}(0^-)\not=1$ and the equation (\ref{bsoleq}) has only one solution : $$b^*=\frac{c(r-\psi(1))}{r\alpha}\frac{1-{\cal{L}}(0^-)}{1-{\cal{G}}(0^-)}=\frac{c(r-\psi(1))}{r\alpha}lim_{x\uparrow 0}\frac{1-{\cal{L}}(x)}{1-{\cal{G}}(x)}.$$
              
               The function $s$  has the form  $g(., B_c)=s(.)$ and is convex, thus it is continuous, in particular it is continuous at $B_c$. We deduce that $B_c=b^*$.
           \end{preuve}
             \hfill  $\Box$
             \\
           
           When $\cal{G}$ is continuous at $x=0$, $B_c$ is more technical to obtain, but it has the same form.
           
            \begin{theoreme}
             \label{th1}
             Let $\sigma>0$. Under Assumptions \ref{h1}, \ref{h2} and \ref{h3_ch1}, we suppose that the function 
             $\cal{G}$ is continuous at $x=0$. 
             \begin{enumerate}
\item
If $\cal{G}$ has  left derivative at $x=0$ (say  ${\cal{G}}'(0^-)$), then  $\cal{L}$ has left derivative at $x=0$ (say  ${\cal{L}}'(0^-)$). 
\item
 If moreover ${\cal{G}}'(0^-)\not=0,$ then $B_c\in [\tilde b , \frac{(r-\psi(1))c}{r\alpha}[$ where $\tilde b=\frac{(r-\psi(1))c}{r\alpha}lim_{x\uparrow 0}\frac{1-{\cal{L}}(x)}{1-{\cal{G}}(x)}$.
 \item
If moreover $g(., \tilde b)$ is strictly convex on $]\tilde b, ~\infty[$,
\end{enumerate}
     then the smaller optimal stopping time is $$\tau^*(c)=inf\{t\geq 0 : V_t\leq 
             B_c\}, \hbox{~where~} B_c=\tilde b.$$
             \end{theoreme}
              \begin{preuve}
             \\
             $(1)$ By Remark \ref{intervalle_ch1},  there exists $B_c$ such that $g(., B_c)=s(.)$. The function $s$ is convex,  therefore the right and  left derivatives exist everywhere and
             \begin{equation}
             \label{deriv_s}
             s'(v^-)\leq s'(v^+) \hbox{~for all~} v\in\R^*_+,
             \end{equation}
             where   $s'(v^-)$ and $s'(v^+)$ are the left and right derivatives of $s$ at $v$. In particular, this means that
             $$g(v,B_c)=\frac{-\alpha v}{r-\psi(1)}{\cal{G}}\left(ln\frac{B_c}{v}\right)+\frac{c}{r}{\cal{L}}\left(ln\frac{B_c}{v}\right)=s(v)$$
             has right and  left derivatives at $v=B_c$.
             Since $\cal{G}$ has right and left derivatives at $x=0$, then $\cal{L}$ has also right and left derivatives at $x=0$.
             \\
             $(2)$ Let us make $v=B_c$ in (\ref{deriv_s}) :
             $$\frac{-\alpha }{r-\psi(1)}\leq \frac{-\alpha }{r-\psi(1)}+\frac{\alpha }{r-\psi(1)}{\cal{G}}'(0^-)-\frac{c}{rB_c}{\cal{L}}'(0^-).$$
                         We deduce that  $B_c\geq\tilde b=\frac{(r-\psi(1))c}{r\alpha}\frac{{\cal{L}}'(0^-)}{{\cal{G}}'(0^-)}=\frac{(r-\psi(1))c}{r\alpha}lim_{x\uparrow 0}\frac{1-{\cal{L}}(x)}{1-{\cal{G}}(x)}$.
   \\
       $(3)$  If moreover    $g(., \tilde b)$ is strictly convex on $]\tilde b, ~\infty[$, then
       \begin{equation}
       \label{gf}
g(v, \tilde b)>f(v) \hbox{~for all~} v> \tilde b.
\end{equation}
       Indeed, the graph of  $f$ is tangent to the graph of $g(.,\tilde b)$ in $v=\tilde b$. 
       
       Suppose that $B_c> \tilde b$, then 
       $$f(B_c)=s(B_c)=g(B_c,B_c)\geq g(B_c, \tilde b)$$    which contradicts    (\ref{gf}).    
              \end{preuve}
             \hfill  $\Box$
             \\
             
             We stress the following consequence of Theorem \ref{th1}.
          \begin{remarque}If ${\cal{G}}'(0^-)$ exists, then ${\cal{L}}'(0^-)$ exists.\end{remarque}
             \begin{proposition}
            For any $\eps>0$, let
             $$\tau_{\eps}(c)=inf\{t\geq 0 : e^{-rt}s(V_t)\leq e^{-rt}f(V_t) 
             +\eps\}.$$
            If  $\sigma>0$, then  under Assumptions \ref{h1}, \ref{h2} and \ref{h3_ch1}, $\pr(\tau_{\eps}(c)<\infty)=1$ and\\
             $lim_{\eps\rightarrow 0}\tau_{\eps}(c)=\tau^*(c)$.
             \end{proposition}
             \begin{preuve}
             \\
             Since $\esp[sup_{t\geq 0}max(e^{-rt}f(V_t),0)]\leq \frac{c}{r}$, then using Lemma ~\ref{tpepsopt} of  Section  \ref{appendix1}, 
             $$\pr(\tau_{\eps}(c)<\infty)=1.$$
The sequence $(\tau_{\eps}(c), ~~\eps\geq 0)$ is a decreasing sequence of stopping times, hence the limit $\tau_0=lim_{\eps\rightarrow 0}\tau_{\eps}(c)$ exists and using Theorem \ref{th5} of Section \ref{appendices}, it is equal to $\tau_0=\tau^*(c)$.
             \end {preuve}
             \hfill  $\Box$
             \\

 \subsubsection {Application to  Finance \\}

             This type of optimal stopping problem can be applied in Finance : the process $V$ describes the assets value of a given firm. The rate  $r$ describes the discount current rate and $\psi(1)=\frac{1}{t}ln \esp\left(\frac{V_t}{V_0}\right)$ the expected asset growth rate. Suppose that the firm generates cash flows  at the rate $\alpha V_t$ at any time $t$. The firm issues bonds and pays coupons indefinitely 
           (meaning that $c$ is a speed of payement). The expected present value of the cash flows generated by the firm until the  liquidation time $\tau$ is 
             $$\esp_v\left[\int_0^{\tau}e^{-rt}(\alpha V_s-c)ds\right].$$
             At a fixed time $t$, the equity owners look for an optimal liquidation time : they want to maximize the expected present value of the cash flows generated by the firm until the  liquidation time $\tau$. This one  corresponds to the solution of the optimal stopping problem (\ref{2}).

             \section {Examples}
\label{section4_ch1}

The examples presented in this section are  some models where  $\esp_v\left(e^{-r\tau_b+aX_{\tau_b}}\right)$ is known for all $b$. Next we consider some particular  Lévy processes, we check that the assumptions of Theorem \ref{th1} or Theorem \ref{th2}  are satisfied  and we solve the problem (\ref{2}) in each case. We start with a continuous Lévy process (Brownian motion), then we continue with a double exponential jump-diffusion process, a particular spectrally positive process and we finish with the  Poisson process.

             \subsection{Brownian motion}
\label{MB}

     We find Duffie and Lando's result, the optimal stopping problem (\ref{2}) being already studied in \cite{Duf} for a Brownian motion with drift. The parameters $(m, \sigma, \delta, (\theta-1)C)$ of Duffie and Lando's model correspond here to $(0, 1, \alpha, c)$. 
     Contrary to their method, our method  allows to avoid  long calculations of the  integro-differential operator.

             Let $X=W$ where $(W_t, ~t\geq 0)$ is a standard Brownian motion. Then $V=ve^W$. In this case $\psi(1)=\frac{1}{2}$ and Assumption  \ref{h1} is checked. We impose (Assumption
             \ref{h2}) that $r>\frac{1}{2}$. 
             \begin{lemme}
             \label{h3_mb}
             The  Brownian motion checks the relation  (\ref{limRn}), i.e.  
             $$lim_{n\rightarrow \infty}\esp\left(e^{-rR_n+W_{R_n}}\1_{R_n<\infty}|W_0=0\right)=0$$
             where  $R_n=inf\{t\geq 0 : e^{-rt+W_t}\geq n\}$. Consequently, Assumption \ref{h3_ch1} is satisfied.
             \end{lemme}

             \begin{preuve} 
             \\
             The process $W$ is continuous, so $$R_n=inf\{t\geq 0 : e^{-rt+W_t}\geq n\}=inf\{t\geq 0 : -rt+W_t=ln(n)\}.$$ 
             Thus $\esp\left(e^{-rR_n+W_{R_n}}\1_{R_n<\infty}|W_0=0\right)=n\pr\left(R_n<\infty|W_0=0\right).$
             
             We apply a result of  \cite{KS} (page 197), and  we obtain 
             $$\esp\left(e^{-rR_n+W_{R_n}}\1_{R_n<\infty}|W_0=0\right)=ne^{-2rln(n)}=n^{1-2r}.$$
             However $r>\frac{1}{2}$, then $lim_{n\rightarrow \infty}\esp\left(e^{-rR_n+W_{R_n}}\1_{R_n<\infty}|W_0=0\right)=0$.
\end {preuve}
\hfill  $\Box$
\\

             The hypothesis of Proposition \ref{formetau} are checked and the smallest optimal stopping time has the form 
             $$\tau_{b_c}=inf \{t\geq 0 : ~V_t \leq b_c\}=inf \{t\geq 0 : ~V_t =b_c\}$$
             because $V$ is a continuous process. 

            Using Remark 8.3 page 96  of \cite{KS} which gives the Laplace transform of a hitting time in the case of a Brownian motion,  $$ {\cal{L}}(x)=e^{x\sqrt{2r}}\hbox{~and~}  {\cal{G}}(x)=e^{x(\sqrt{2r}+1)}, ~x<0.$$
The function $g$ has the form 
              \begin {equation}
             \nonumber
             g(v,b)=\left\{
             \begin {array}{l@{\quad}l}
             -\frac{\alpha v}{r-\frac{1}{2}}+ \frac{c}{r}& if \ v\leq b\\
             \left(\frac{-\alpha 
             b}{r-\frac{1}{2}}+\frac{c}{r}\right)\left(\frac{v}{b}\right)^{-\sqrt{2r}} & if \ v>b. \\
             \end{array}
             \right.
             \end{equation}
             
             We notice that the function $x\mapsto {\cal{G}}(x)$ is continuous at $x=0$.
Let us check the hypothesis of Theorem \ref{th1} :
             \\
             \\
             (1) $\cal{G}$ has  left derivative at $x=0$.
\\
(2) Moreover ${\cal{G}}'(0^-)=\sqrt{2r}+1\not=0.$ 
Since ${\cal{L}}'(0^-)=\sqrt{2r}$ then  $\tilde b=\frac{c\sqrt{2r}(r-\frac{1}{2})}{\alpha  r(\sqrt{2r}+1)}$.
\\
(3)  Remark that   the function $g(., \tilde b)$ belongs to  
             $\cal{C}$$^2(]\tilde b, ~\infty[)$. Its second derivative is equal to  
             $$\frac{\partial^2 g}{\partial v^2}(v, \tilde b)=\left(\frac{-\alpha \tilde b}{r-\frac{1}{2}}+
          \frac{c}{r}\right)\frac{\sqrt{2r}(\sqrt{2r}+1)}{\tilde b^2}\left(\frac{v}{\tilde b}\right)^{-\sqrt{2r}-2} =\frac{c\sqrt{2r}}{r\tilde b^2}\left(\frac{v}{\tilde b}\right)^{-\sqrt{2r}-2}>0$$ 
             and the function  $g(., \tilde b)$ is strictly convex on $]\tilde b, ~\infty[$.

             Thus $B_c=\tilde b$ and the optimal stopping time is $$\tau^*(c)=inf \{t\geq 0 : ~V_t =B_c\}.$$

             \begin{proposition}
             If $X$ is a standard Brownian motion, then with the notations introduced in Section  \ref{MB},
             \begin{enumerate}
             \item
             The smallest optimal stopping time is $\tau^*(c)=inf \{t\geq 0 : ~V_t =B_c\}$ where \\ $B_c=\frac{c\sqrt{2r}(r-\frac{1}{2})}{\alpha r(\sqrt{2r}+1)}$.
             \item
       The value function  $w$ is given by $$w(v)=\frac{\alpha 
             v}{r-\frac{1}{2}}-\frac{c}{r}\left(1+\frac{1}{\sqrt{2r}}\left(\frac{v}{B_c}\right)^{-\sqrt{2r}}\right) \hbox{~~where~} v>B_c.$$
             \end{enumerate}
             \end{proposition}

             \subsection{Double exponential jump-diffusion process}
\label{KW}

Using Lemma \ref{lem3}, the problem (\ref{2}) can be brought back to an optimal stopping problem for an American Put option with strike price  $\frac{c(r-\psi(1))}{r\alpha}$ :
$$s(v)=sup_{\tau\geq 0}\esp_v\left[e^{-r\tau} \left(\frac{-\alpha  V_{\tau}}{r-\psi(1)} +\frac{c}{r}\right)^+\right] =\frac{\alpha}{r-\psi(1)}sup_{\tau\geq 0}\esp_v\left[e^{-r\tau} \left(-V_{\tau}+\frac{c(r-\psi(1))}{r\alpha}\right)^+\right].$$
We find the result of Theorem 1 of \cite{Kou2}. In  \cite{Kou2}, the authors solve an optimal stopping problem for an American Put option. They use the Wiener-Hopf factorization. But, in general explicit calculation of the Wiener-Hopf factorization is difficult. Because of the memoryless property of the exponential distribution, they can solve the problem explicitly. Our method is much  easier to use than their method,  much more rapid and it can be used for any Lévy process.

In \cite{pham}, Pham  studies an optimal stopping problem for an American Put option with finite time horizon. His model is a jump-diffusion one and the jumps are not restricted to any particular law. He uses integro-differential equations to solve his problem.

            The next model is   Kou and Wang's model (we refer to \cite {Kou, Kou2} and \cite{Chen}). Indeed, we suppose that $X$ is a mixed diffusion-jump process and the jump size is a double exponential distributed random variable   :
              \begin {equation}
             \label{double_exp_ch1}
             X_t=mt+\sigma W_t+\sum_{i=1}^{N_t} Y_i, ~t\geq 0,
             \end {equation}
             where  $(W_t, ~t\geq 0)$ is a standard Brownian motion, $(N_t, ~t\geq 0)$ is  a Poisson process with constant positive intensity $a$, $(Y_i, ~i\in \N )$ is a sequence of independent and identically distributed random variables. The common density of $Y$ is given by
             $$f_Y(y)=p\eta_1e^{-\eta_1y}\1_{y>0}+q\eta_2e^{\eta_2y}\1_{y<0}, ~~y\in\R,$$
             where $p+q=1$, $p,q>0$, $\eta_1>1$ and $\eta_2>0$.
            Moreover we suppose that $(Y_i, ~i\in \N )$, $(N_t, ~t\geq 0)$ and $(W_t, ~
             t\geq 0)$ are independent.
             
             We treat separately the cases where $q=0$ or $p=0$.

            The condition  $\eta_1>1$ implies $\esp(e^Y)<\infty$ and $\esp_v(V_t)<\infty$ for every $t\geq 0$. Assumption \ref{h1} is checked  and $\psi(1)=m+\frac{\sigma^2}{2}+a\esp(e^Y-1)$. Here Assumption \ref{h2} is \\ $r>m+\frac{\sigma^2}{2}+a\esp(e^Y-1)$ where $\esp(e^Y-1)=\frac{\eta_1p}{\eta_1-1}+\frac{\eta_2q}{\eta_2+1}-1$.
             
             \begin{lemme}
             \label{h3_doubleexp}
             The process $X$ introduced in (\ref{double_exp_ch1})  checks the relation  (\ref{limRn}), i.e. 
             $$lim_{n\rightarrow \infty}\esp\left(e^{-rR_n+X_{R_n}}\1_{R_n<\infty}|X_0=0\right)=0$$
            where $R_n=inf\{t\geq 0 : e^{-rt+X_t}\geq n\}$. Consequently, Assumption \ref{h3_ch1} is satisfied.
             \end{lemme}
             
            The demonstration of this lemma rests on Corollary 3.3  of \cite{Kou} recalled in Section  \ref{appendix2} (we refer to  relation (\ref{1_ch4})).

\begin{preuve}  
\\
Remark that $R_n=inf\{t\geq 0 : (-r+m)t+\sigma W_t+\sum_{i=1}^{N_t} Y_i\geq ln(n)\}$, $n\in\N^*$. 
\\
We apply (\ref{1_ch4}) of Section \ref{appendix2} to $r=0$, $\beta=1$,  $b=ln(n)$ and \\ $X_t=(-r+m)t+\sigma W_t+\sum_{i=1}^{N_t} Y_i$ (in fact we replace the drift  $m$ by $-r+m$) :
\begin{equation}
\label{lemme1.3.4_1_ch1}
\esp\left(e^{-rR_n+X_{R_n}}\1_{R_n<\infty}|X_0=0\right)
=n^{1-\psi_1}\frac{(\eta_1-\psi_1)(\psi_0-1)}{(\psi_0-\psi_1)(\eta_1-1)}+n^{1-\psi_0}\frac{(\psi_0-\eta_1)(\psi_1-1)}{(\psi_0-\psi_1)(\eta_1-1)}
\end{equation}
where $0<\psi_1<\eta_1<\psi_0<\infty$ are the positive roots of equation $f(\psi)=0$ with
             \begin{equation}
             \nonumber      f(\psi)=(-r+m)\psi+\frac{\sigma^2}{2}\psi^2+a\left[\frac{\eta_1p}{\eta_1-\psi}+\frac{\eta_2q}{\eta_2+\psi}-1\right].
             \end{equation}
             
  The equation $f(\psi)=0$ has exactly four roots (see Lemma 2.1 of \cite{Kou}), and the two positive solutions are both strictly greater than 1. Indeed,         
             \begin{center}
             \begin{tabular}{|l|c r c l c r|}
   \hline
   $\psi$ & $-\infty$ & $-\eta_2$ & $0$ & $1$ & $\eta_1$ & $\infty$ \\
   \hline
   $f(\psi)$ & $\infty$ & $-\infty$ \vline $\infty$ & $0$ & $f(1)$ & $\infty$ \vline $-\infty$ & $\infty$ \\
   \hline
\end{tabular}
\end{center}
where $f(1)=(-r+m)+\frac{\sigma^2}{2}+a\left[\frac{\eta_1p}{\eta_1-1}+\frac{\eta_2q}{\eta_2+1}-1\right]~=-r+m+\frac{\sigma^2}{2}+aE(e^Y-1)<0$ according to Assumption \ref{h2} rewritten in this case.

Since $1<\psi_1<\eta_1<\psi_0$, then by taking the limit in (\ref{lemme1.3.4_1_ch1}) we obtain : $$lim_{n\rightarrow \infty}\esp\left(e^{-rR_n+X_{R_n}}1_{R_n<\infty}|X_0=0\right)=0.$$
\end {preuve}
\hfill  $\Box$
\\     

            The assumptions of Proposition \ref{formetau} are checked and the smallest optimal stopping time has the form  $$\tau_{b_c}=inf \{t\geq 0 : ~V_t \leq b_c\}.$$ 
            The results of \cite{Kou} allow us to write for all $x<0$  the form of the functions $x\mapsto {\cal{L}}(x)$, $x\mapsto{\cal{G}}(x)$ and therefore $g(., b)$ on  $]b,~\infty[$ :
 \begin{align*}           {\cal{L}}(x)=&\frac{\psi_2(\eta_2+\psi_3)}{(\psi_2-\psi_3)\eta_2}e^{-x\psi_3}-\frac{\psi_3(\eta_2+\psi_2)}{(\psi_2-\psi_3)\eta_2}e^{-x\psi_2},\\
{\cal{G}}(x)=&e^{x}\left[ \frac 
             {(\eta_2+\psi_3)(\psi_2-1)}{(\psi_2-\psi_3)(\eta_2+1)}e^{-x\psi_3}+ 
             \frac{(\eta_2+\psi_2)(1-\psi_3)}{(\psi_2-\psi_3)(\eta_2+1)}e^{-x\psi_2}\right],
             \end{align*}
where  $-\infty<\psi_3<-\eta_2<\psi_2<0$ are the two negative roots  of the equation
             \begin{align*}
             m\psi+\frac{\sigma^2}{2}\psi^2+a[\frac{\eta_1p}{\eta_1-\psi}+\frac{\eta_2q}{\eta_2+\psi}-1]=r.
             \end{align*}

             Thus, the function  $g$ is equal to :
             
             \begin {equation}
             \nonumber
             g(v,b)=\left\{
             \begin {array}{l@{\quad}l}
             -\frac{\alpha v}{r-\psi(1)}+\frac{c}{r}& if \ v\leq b\\
  A_{b}(\frac{v}{b})^{\psi_3}+D_{b}(\frac{v}{b})^{\psi_2} 
             & if \ v>b,\\
             \end{array}
             \right.
             \end{equation}
             \begin{align*}
              \hbox{where } ~~~~A_{b}=&-\frac{\alpha 
             b}{r-\psi(1)}\frac{(\eta_2+\psi_3)(\psi_2-1)}{(\psi_2-\psi_3)(\eta_2+1)}+\frac {c}{r}  
             \frac{\psi_2(\eta_2+\psi_3)}{(\psi_2-\psi_3)\eta_2} ~~~~\hbox{ and} 
              \\
              D_{b}=&-\frac{\alpha 
             b}{r-\psi(1)}\frac{(\eta_2+\psi_2)(1-\psi_3)}{(\psi_2-\psi_3)(\eta_2+1)}-\frac{c}{r}\frac{\psi_3(\eta_2+\psi_2)}{(\psi_2-\psi_3)\eta_2}.
              \end{align*}
              
             Remark that the function $x\mapsto {\cal{G}}(x)$ is continuous at $x=0$.
 Let us check the assumptions of Theorem \ref{th1} :
\\
             (1) $\cal{G}$ has  left derivative at $x=0$.
\\
(2) Moreover ${\cal{G}}'(0^-)=\frac{(1-\psi_2)(1-\psi_3)}{(\eta_2+1)}\not=0.$ 
Since ${\cal{L}}'(0^-)=\frac{\psi_2\psi_3}{\eta_2}$ then  $\tilde b=\frac{c(r-\psi(1))\psi_2\psi_3(\eta_2+1)}{r\alpha\eta_2(1-\psi_2)(1-\psi_3)}.$
\\
(3) Remark that $g(., \tilde b)\in \cal{C}$$^2(]\tilde b, ~\infty[)$ and 
             \begin{align*} 
             \frac{\partial^2 g}{\partial 
             v^2}(v, \tilde b)=
           &A_c(\frac{1}{\tilde b})^{\psi_3}\psi_3(\psi_3-1)v^{\psi_3-2}+
           D_c(\frac{1}{\tilde b})^{\psi_2}\psi_2(\psi_2-1)v^{\psi_2-2}\\
             \hbox{where } 
             ~~~~A_c=&\frac{c\psi_2(\eta_2+\psi_3)}{r(\psi_2-\psi_3)\eta_2(1-\psi_3)}>0,\\ 
             D_c=&-\frac{c\psi_3(\eta_2+\psi_2)}{r(\psi_2-\psi_3)\eta_2(1-\psi_2)}>0.
             \end{align*}
             Thus $\frac{\partial^2 g}{\partial v^2}(v,\tilde b)>0$ for $v>\tilde b$ and $g(., \tilde b)$ is strictly convex on $]\tilde b, ~\infty[$.

             We can apply Theorem \ref{th1}, $B_c=\tilde b$ and the optimal stopping time is 
             $$\tau^*(c)=inf\{t\geq 0 : V_t\leq B_c\}.$$
             
             \begin{proposition}
             Let $X$  be the process introduced in (\ref{double_exp_ch1}). Then, using the notations introduced in Section  \ref{KW} :
             \begin{enumerate}
             \item
              The smallest optimal stopping time is $\tau^*(c)=inf \{t\geq 0 : ~V_t \leq B_c\}$ where \\$B_c=\frac{c(r-\psi(1))\psi_2\psi_3(\eta_2+1)}{r\alpha\eta_2(1-\psi_2)(1-\psi_3)}$.
             \item
               For $v>B_c$, the value function $w$ is equal to
             \begin {equation}
             \nonumber
            w(v)= \frac{\alpha  v}{r-\frac{1}{2}}-\frac{c}{r}+ \frac{c\psi_2(\eta_2+\psi_3)}{r(\psi_2-\psi_3)\eta_2(1-\psi_3)}\left(\frac{v}{B_c}\right)^{\psi_3}-\frac{c\psi_3(\eta_2+\psi_2)}{r(\psi_2-\psi_3)\eta_2(1-\psi_2)}\left(\frac{v}{B_c}\right)^{\psi_2}.
             \end {equation}
             \end{enumerate}

             \end{proposition}

 \subsection{Exponential jump-diffusion process}
\label{sautspos}

To our knowledge, the following result seems to be new.

            In this section, we suppose that $X$ is a mixed diffusion-jump process and the jump size is a random variable with an  exponential distribution :
              \begin {equation}
             \label{exp_ch1}
             X_t=mt+\sigma W_t+\sum_{i=1}^{N_t} Y_i, ~t\geq 0,
             \end {equation}
             where  $(W_t, ~t\geq 0)$ is a standard Brownian motion, $(N_t, ~t\geq 0)$ is a Poisson process with constant positive intensity $a$, $(Y_i, ~i\in \N )$ is a sequence of independent and identically distributed random variables with an exponential distribution, i.e. the common density of $Y$ is given by $f_Y(y)=\eta_1e^{-\eta_1y}\1_{y>0}$ where
 $\eta_1>1$. Moreover we suppose that $(Y_i, ~i\in \N )$, $(N_t, ~t\geq 0)$ and $(W_t, ~t\geq 0)$ are independent. This is a particular Lévy process with positive jumps.
 
            As for the double exponential jump-diffusion process (here we consider $q=0$, so $p=1$),  the condition $\eta_1>1$ implies $\esp(e^Y)<\infty$.  Assumption \ref{h1} is thus checked. In this particular case, $\psi(1)=m+\frac{\sigma^2}{2}+\frac{a}{\eta_1-1}$ and Assumption \ref{h2} is $r>m+\frac{\sigma^2}{2}+\frac{a}{\eta_1-1}$.
             
             \begin{lemme}
             \label{h3_exp}
             The process $X$ introduced in (\ref{exp_ch1}) satisfies the relation (\ref{limRn}), i.e. 
             $$lim_{n\rightarrow \infty}\esp\left(e^{-rR_n+X_{R_n}}\1_{R_n<\infty}|X_0=0\right)=0$$
            where $R_n=inf\{t\geq 0 : e^{-rt+X_t}\geq n\}$. Consequently, Assumption \ref{h3_ch1} is checked.
             \end{lemme}
             
            We use   (\ref{1_ch4}) and Remark \ref{p1} of Section \ref{appendix2} to prove this result. 
           The demonstration of this lemma is the same as that of Lemma \ref{h3_doubleexp}, the only difference is that here  $p=1$ (thus $q=0$) and $0<\psi_1<\eta_1<\psi_0<\infty$ are the positive roots of equation $f(\psi)=0$ where
            $ f(\psi)=(-r+m)\psi+\frac{\sigma^2}{2}\psi^2+\frac{a\psi}{\eta_1-\psi}.$

            The assumptions of Proposition \ref{formetau} are checked and the smallest optimal stopping time has the form  $$\tau_{b_c}=inf \{t\geq 0 : ~V_t \leq b_c\}.$$ 
            
            The following result of \cite{Em} allows us to write for every $x<0$  the form of the functions $x\mapsto{\cal{L}}(x)$ and $x\mapsto{\cal{G}}(x)$  :
            
              \begin{proposition} (Theorem 1 of \cite{Em})
             \label{emery}

            Let $X$ be a Lévy process with no negative jumps and characteristic function $\Psi$ :\\
             $\esp(e^{i\lambda X_t})=e^{-t\Psi(\lambda)}.$  
      \begin{enumerate}     
\item
The equation $r+\Psi(\lambda)=0$ has at most one root in $Im(\lambda)>0$.
\item Let $x<0$ and $\tau=inf\{t\geq 0 : X_t\leq x\}$. If $r+\Psi(\lambda)=0$ has a root denoted $i\bar{\lambda}$, $\bar{\lambda}>0$, then
             $$\esp\left(e^{-r\tau+iq X_{\tau}}\right)=e^{x\bar{\lambda}+ixq}.$$
             \end{enumerate}
             \end{proposition}

In our particular case $\Psi(\lambda)=-im\lambda+\frac{\sigma^2\lambda^2}{2}-\frac{ai\lambda}{\eta_1-i\lambda}$. Since the equation $r+\Psi(\lambda)=0$ has a root in $Im(\lambda)>0$, then we can apply Proposition
             \ref{emery}   :
             $${\cal{L}}(x)=e^{x\bar{\lambda}} \hbox{~and ~} {\cal{G}}(x)=e^{x(\bar{\lambda}+1)}. $$
    
    The function $g$ has the form 
        \begin {equation}
             \nonumber
             g(v,b)=\left\{
             \begin {array}{l@{\quad}l}
             -\frac{\alpha v}{r-\psi(1)}+ \frac{c}{r}& if \ v\leq b\\
             \left(\frac{-\alpha 
             b}{r-\psi(1)}+\frac{c}{r}\right)\left(\frac{b}{v}\right)^{\bar{\lambda}} & if \ v>b. \\
             \end{array}
             \right.
             \end{equation}
             
             Remark that the function $x\mapsto {\cal{G}}(x)$ is continuous at $x=0$.
Let us check the assumptions of Theorem \ref{th1} :
          \\
             (1) $\cal{G}$ has  left derivative at $x=0$.
\\
(2) Moreover ${\cal{G}}'(0^-)=\bar{\lambda}+1\not=0.$ 
Since ${\cal{L}}'(0^-)=\bar{\lambda}$ then  $\tilde b=\frac{c\bar{\lambda}(r-\psi(1))}{\alpha 
             r(\bar{\lambda}+1)}.$
\\   
   (3) Remark  that $g(., \tilde b)\in \cal{C}$$^2(]\tilde b, ~\infty[)$ and 
             $\frac{\partial^2 g}{\partial 
             v^2}(v, \tilde b)=\frac{c\bar{\lambda}}{r\tilde b^2}
             \left(\frac{\tilde b}{v}\right)^{\bar{\lambda}+2}>0$. Thus the function $g(., \tilde b)$ is strictly convex on $]\tilde b, ~\infty[$.

             Then $B_c=\tilde b$. Since, $X$ has no negative jumps, the smallest optimal stopping time is  $$\tau^*(c)=inf \{t\geq 0 : ~V_t =B_c\}.$$

             \begin{proposition}
            Let $X$ be the process introduced in (\ref{exp_ch1}). Then, using the notations introduced in section \ref{sautspos}, we have :
             \begin{enumerate}
             \item
              The smallest optimal stopping time is $\tau^*(c)=inf \{t\geq 0 : ~V_t \leq B_c\}$ where \\ $B_c=\frac{c\bar{\lambda}(r-\psi(1))}{\alpha r(\bar{\lambda}+1)}$.
              \item
             For $v>B_c$, the value function $w$ is equal to 
             $w(v)=\frac{\alpha 
             v}{r-\psi(1)}-\frac{c}{r}+\frac{c}{r(1+\bar{\lambda})}\left(\frac{B_c}{v}\right)^{\bar{\lambda}}.$
             \end{enumerate}
 \end{proposition}

             \subsection{ Poisson process}
\label{poisson}

In absence of a known reference, we treat the easy case $X=-N$ where $(N_t, ~t\geq 0)$ is a Poisson process with constant positive intensity $a$. Then $V=ve^{-N}$. In this case the Laplace transform of the process $X$ exists, Assumption \ref{h1} is checked and
             $\psi(1)=a(e^{-1}-1)$. We note that here Assumption  \ref{h2} is always checked since  $r>0>a(e^{-1}-1).$
             Moreover, the process $t\mapsto Y_t=e^{-rt}\left(\frac{-\alpha ve^{-N_t}}{r-a(e^{-1}-1)}+\frac{c}{r}\right)$ is bounded, thus  of class $D$. 
             
             In this particular case, since the Gaussian component of  $X$ is null,
Lemma \ref{lem3} cannot be applied. Let us recall  that this lemma allows us to rewrite the function $s$ in a particular form ($s(v)=sup_{\tau\geq 0}\esp_v[e^{-rt}f^+(V_t)]$) in order  to be able to apply Theorem \ref{th5} of Section \ref{appendices} and to find  the smallest optimal stopping time form. However, since the process $t\mapsto Y_t=e^{-rt}f(V_t)$ is bounded, then  Theorem \ref{th5} of Section \ref{appendices} may be applied directly. The smallest optimal stopping time is $\tau^*=inf\{t : f(V_t)=s(V_t)\}$. In this particular case where $X = -N$, the function $s$ is defined  on $\{v, ve^{-1}, ve^{-2}...\}$, its continuous prolongation by linear interpolation is convex and the conclusion of Proposition \ref{formetau} is true.

The smallest optimal stopping time has the form
$$\tau_{b_c}=inf \{t\geq 0 : ~V_t \leq b_c\}=inf \left\{t\geq 0 : ~N_t \geq 
             ln\left(\frac{v}{b_c}\right)\right\}$$ and in this case it coincides with  a jump time of the process $N$.
             
             Let $(T_i, ~i\in\N^*)$ be a sequence of random variables describing the jump times of the process  $N$.

             To find the form of the functions ${\cal{L}}$ and ${\cal{G}}$ we calculate
             $\esp\left(e^{-r\tau_x}\right)$ and 
             $\esp\left(e^{-r\tau_x-N_{\tau_x}}\right)$ where $\tau_x=inf \{t\geq 0 : ~N_t \geq -x\}$:
             $${\cal{L}}(x)=\sum_{i\geq 0}\esp\left(e^{-rT_i}\1_{\tau_x=T_i}\right)=\1_{x\geq 0}+\sum_{i\geq 
             1}\esp\left(e^{-rT_i}\right)\1_{i-1<-x\leq i}.$$
             Since $T_i=S_1+...+S_i$ where $(S_j, ~j\in\N)$ is a sequence of independent and identically distributed random variables with  exponential distribution with parameter $a$, then 
             $\esp\left(e^{-rT_i}\right)=\esp\left(e^{-rS_1}\right)^i=\left(\frac{a}{r+a}\right)^i$ and 
             $${\cal{L}}(x)=\1_{x\geq 0}+\sum_{i\geq 
             1}\left(\frac{a}{r+a}\right)^i\1_{i-1<-x\leq i}.$$

             In the same way, we calculate  
            $${\cal{G}}(x)=\sum_{i\geq 
             0}\esp\left(e^{-rT_i-i}\1_{\tau_x=T_i}\right)=\1_{x\geq 0}+\sum_{i\geq 
             1}\left(\frac{a}{e(r+a)}\right)^i\1_{i-1<-x\leq i}.$$

             The function $g(., b)$ has the form :
             \begin {equation}
             \nonumber
             g(v,b)=\left\{
             \begin {array}{l@{\quad}l}
             -\frac{\alpha v}{r-\psi(1)}+\frac{c}{r}& if \ v\leq b\\
             -\frac{\alpha v}{r-\psi(1)}\frac{a}{e(r+a)}+\frac{c}{r}\frac{a}{r+a} 
             & if \ b<v\leq be\\
             -\frac{\alpha 
             v}{r-\psi(1)}\left(\frac{a}{e(r+a)}\right)^2+\frac{c}{r}\left(\frac{a}{r+a}\right)^2 & if \ be<v\leq be^2\\
             ...
             \end{array}
             \right.
             \end{equation}

             Remark that ${\cal{G}}$ is discontinuous at $x=0$. By Theorem \ref{th2}, $B_c=\frac{ce(r-\psi(1))}{\alpha(er+ea-a)}$ and  the optimal stopping time is 
             $$\tau^*(c)=inf \{t\geq 0 : ~V_t \leq B_c\}.$$

             \begin{proposition}
             Let $X=-N$ where $(N_t, ~t\geq 0)$ is a Poisson process with constant positive intensity $a$. Then with the notations introduced in Section \ref{poisson} : 
              \begin{enumerate}
              \item
              The smallest optimal stopping time is $\tau^*(c)=inf \{t\geq 0 : ~V_t \leq B_c\}$ where \\ $B_c=\frac{ce(r-\psi(1))}{\alpha(er+ea-a)}$.
              \item
              For $v>B_c$, the value function $w$ is equal to  $$w(v)=\sum_{i\geq 1}\left[\frac{\alpha 
           v}{r-\psi(1)}\left(1-\left(\frac{a}{e(r+a)}\right)^i\right)-\frac{c}{r}\left(1-\left(\frac{a}{r+a}\right)^i\right)\right]\1_{B_ce^{i-1}<v\leq 
             B_ce^i}.$$
             \end{enumerate}
             \end{proposition}

             \section{Spectrally negative processes. Case of an unbounded variation process.}
               \label{levneg}
              
              The case of spectrally negative processes is more complicated and requires more calculations than the other examples. We present here the example of a  spectrally negative process  with a non null Gaussian component.

             \subsection{Notations, hypothesis and tools}
             \label{sautneg_not}
            Throughout this section, we suppose that  $X$ is a real-valued Lévy process with no positive jumps. Some authors say that $X$ is spectrally negative. The case when $X$ is either a negative Lévy process with decreasing paths or a deterministic drift are excluded in this sequel. The Laplace transform of such a process exists and has the following  form (see \cite{Bert1} page 187-189 for details) :
         
             \begin{equation}
             \nonumber
             \esp(e^{\lambda X_t})=e^{t\psi(\lambda)} 
             ~~~\hbox{where}~~~\lambda \in\R_+ , ~\psi(\lambda)=m\lambda+\frac{\sigma^2}{2}\lambda^2+\int_{-\infty}^0(e^{\lambda 
             x}-1-\lambda x\1_{x>-1})\Pi(dx).
             \end{equation}
             
             The function $\psi : [0, \infty[\rightarrow \R$ is strictly convex and $lim_{\lambda \rightarrow\infty}\psi(\lambda)=\infty.$ We denote by $\Phi(0)$ the largest solution of the equation  $\psi(\lambda)=0$. Observe that $0$ is always a solution of $\psi(\lambda)=0$. If
             $\Phi(0)>0$, then by strict convexity, $0$ and
             $\Phi(0)$ are the only solutions of $\psi(\lambda)=0$. In all cases,
             $\psi :[\Phi(0), \infty[\rightarrow \R_+$ is continuous and increasing, it is a bijection and its inverse is $\Phi :[0, \infty[\rightarrow[\Phi(0), \infty[$ : $\psi\circ \Phi(\lambda)=\lambda$ $(\lambda>0)$. Assumption \ref{h2}, i.e.  $r>\psi(1)$ is equivalent to $\Phi(r)>1$.
             

In the sequel we suppose that $r>\psi(1)$ (Assumption \ref{h2}).

             The process $t\mapsto e^{cX_t-\psi(c)t}$ is a martingale. Define for each $c\geq 0$ the change of measure : 
             $\frac{d\pr^c}{d\pr}\mid_{\F_t}=e^{cX_t-\psi(c)t}$. Following \cite{Kyp}, we introduce a new function  :
             \begin{Defi}
             \label{def0_ch1}
             For any $q\geq 0$ let $W^{(q)} : \R\rightarrow 
             \R_+$   be the function defined by  $$W^{(q)}(x)=e^{\Phi(q)x}\pr^{\Phi(q)}(inf_{t\geq 0}X_t \geq 0|X_0=x).$$
             \end{Defi}

The following proposition presents some properties of the function $W^{(q)}$ shown in  \cite{Kyp} and \cite{CK} :

             \begin{proposition}
             \label{pptes_w}
             For any $q\geq 0$, the function $W^{(q)}$ has the following properties :
              \begin{enumerate}
\item
            $W^{(q)}(x)=0$ for $x<0$ and $W^{(q)}$ is a strictly increasing and continuous function on $[0, \infty[$ whose Laplace transform satisfies 
 $$\int_0^{\infty}e^{-\beta x}  W^{(q)}(x)dx=\frac{1}{\psi(\beta)-q}$$ for  $\beta>\Phi(q)$ (Theorem 8.1 page 214 of
             \cite{Kyp}).
             \item
             $W^{(q)}(0)=0$ if and only if $X$ has unbounded variation (Lemma 8.6 page 223 of  \cite{Kyp}).\item
             $W^{(q)'}(0)=\frac{2}{\sigma^2}$ if   $X$ has unbounded variation (Exercice 8.5 page 234 of \cite{Kyp}).\item
            If the process $X$ has a non null Gaussian component, then the function $W^{(q)}$ belongs to $\cal{C}$$^2 (]0, \infty[)$ (Theorem 2 of
             \cite{CK}). 
          \end{enumerate}    
         \end{proposition}

             We introduce the following functions:
             \begin{Defi}
             \label{def1}
             \begin{enumerate}
             \item
             For any $q\geq 0$ let $Z^{(q)} : \R_+ \rightarrow \R_+$  be the function defined by  $$Z^{(q)}(x)=1+q\int_0^xW^{(q)}(y)dy.$$
            \item
              For any $q\in\C$ and $c\in \R$ such that $\psi(c)\leq q$, let 
             \begin{equation}
             \nonumber
             W^{(q-\psi(c))}_c(x)=e^{-cx}W^{(q)}(x) 
             ~~~~\hbox{and}~~~~Z^{(q-\psi(c))}_c(x)=1+(q-\psi(c))\int_0^xW^{(q-\psi(c))}_c(y)dy
             \end{equation}
             for every $x\geq 0$.
             \end{enumerate}
             \end{Defi}

          \begin{remarque}
          \label{valen0}
          i) From Proposition \ref{pptes_w} (2) and Definition \ref{def1} (2), if $X$ has unbounded variation, then $$W^{(q-\psi(c))}_c(0)=0.$$ 
          ii) Moreover, Definition \ref{def1} (1) and (2) taken in  $x=0$ give $$Z^{(q)}(0)=1, ~Z^{(q-\psi(c))}_c(0)=1.$$ 
          \end{remarque}

             \subsection{Optimal stopping time}

             Let us suppose that $\sigma>0$ ; one can express the Lévy process as the sum of a Brownian motion with drift and a pure  jump process with negative jumps ; there are thus processes with unbounded variation. The case $p=0$ of Section \ref{KW} is a particular case of a negative Lévy process with a non null Gaussian component.
             
             At our knowledge, there does not exist any proof to show that under Assumption \ref{h2}, the process $\left(e^{-rt+X_t}, ~t\geq 0\right)$ is of class D.
\begin{lemme}
             \label{h3_sautneg}
            Let $X$ be a  Lévy process with a non null Gaussian component satisfying  Assumption \ref{h2}. This process checks the relation (\ref{limRn}), i.e.
             $$lim_{n\rightarrow \infty}\esp\left(e^{-rR_n+X_{R_n}}\1_{R_n<\infty}|X_0=0\right)=0$$
             where  $R_n=inf\{t\geq 0 : e^{-rt+X_t}\geq n\}$.  Consequently, Assumption \ref{h3_ch1} is satisfied.
             \end{lemme}
             
            The proof of this lemma rests on the following result  :
             \begin{lemme}
             \label{kyp_sautneg}
            Let $X$ be a spectrally negative process and for each $x>0$ let $\tau^x$ be the following stopping time $\tau^x=inf\{t>0 : X_t>x\}$. Then 
            \begin{enumerate}
            \item
           For any $x>0$,  $\tau^x=inf\{t>0 : X_t\geq x\}$ $~~\pr(.|X_0=0)$-almost surely (Lemma 49.6 page 373 of \cite{sato}).
            \item $\pr\left(\tau^x<\infty|X_0=0\right)=e^{-\Phi(0)x}$ where $\Phi(0)$ is introduced in Section \ref{sautneg_not}  (Corollary 3.13 page 82 of \cite{Kyp}).
            \item
           Since the process $X$ has no positive jumps, then $\pr(X_{\tau^x}=x|\tau^x<\infty, ~X_0=0)=1$ (page 212 of  \cite{Kyp}).
           \end{enumerate}
             \end{lemme}
             \begin{preuve} of Lemma \ref{h3_sautneg}
             \\
 The stopping time $R_n$ can be written as $R_n=inf\{t\geq 0 : -rt+X_t\geq ln(n)\}$. Since $n>0$, then $\pr(R_n=0|X_0=0)=0$ and $R_n=inf\{t>0 : -rt+X_t\geq ln(n)\}$.
            We apply  Lemma \ref{kyp_sautneg} to $t\mapsto -rt+X_t$ and $x=ln(n)$ :  $R_n=\tau^{ln(n)}$  $~~\pr(.|X_0=0)$-almost surely and 
            \begin{equation}
             \label{lemme1.4.3_1_ch1}
             \esp\left(e^{-rR_n+X_{R_n}}\1_{R_n<\infty}|X_0=0\right)=ne^{-\bar{\Phi}(0)ln(n)}=n^{1-\bar{\Phi}(0)},
             \end{equation}
 where $\bar{\Phi}(0)$ is the largest solution of the equation $\bar{\psi}(\lambda)=0$ where
         $$ \bar{\psi}(\lambda)=(m-r)\lambda+\frac{\sigma^2}{2}\lambda^2+\int_{-\infty}^0(e^{\lambda 
             x}-1-\lambda x\1_{x>-1})\Pi(dx).$$
             
            The largest solution of the equation $\bar{\psi}(\lambda)=0$ satisfies $\bar{\Phi}(0)>1$. Indeed, $\bar{\psi}$ is a continuous function,   $\bar{\psi}(1)=m-r+\frac{\sigma^2}{2}+\int_{-\infty}^0(e^x-1-x\1_{x>-1})\Pi(dx)<0$ (by Assumption  \ref{h2}) and this function goes to infinity when  $\lambda$ goes to infinity.
             
             Let us take the limit in (\ref{lemme1.4.3_1_ch1}),  $$lim_{n\rightarrow \infty}\esp\left(e^{-rR_n+X_{R_n}}\1_{R_n<\infty}|X_0=0\right)=0.$$  
             \end {preuve}
\hfill  $\Box$
\\  

The assumptions of Proposition \ref{formetau} are checked and the smallest optimal stopping time has the form
         $$\tau_{b_c}=inf \{t\geq 0 : ~V_t \leq b_c\}.$$
         
         Now we seek ${\cal{L}}$ and ${\cal{G}}$. In \cite{Kyp}, the author calculates the Laplace transform of a stopping time of the form $\tau^-_x=inf \{t> 0 :  ~X_t < x\}$ and the joint Laplace transform  of $(\tau^-_x, ~X_{\tau^-_x})$. In order to be able to apply his results, we show that if there exists a non null Gaussian component, then the hitting time $\tau^-_x$ is equal to $\tau_x=inf \{t\geq 0 : ~X_t \leq x\}$.
         
             \begin{proposition}
             \label{egalite_tau}
             Let $X$ be a spectrally negative process with a non null Gaussian component.  Then, $\pr(.|X_0=0)$-almost surely
             $$inf \{t\geq 0 : ~X_t \leq x\}=inf \{t\geq 0 : ~X_t < x\}=inf \{t> 0 : ~X_t < x\}$$ for $x\leq 0$.
             \end{proposition}
             
             The proof of this proposition rests on the following result :
             
             \begin{proposition} 
             \label{trlap}
            Let $X$ be a spectrally negative process and $\tau^-_y=inf\{t> 0 : X_t<y\}$. 
            \begin{enumerate}
            \item (Theorem 8.1 (ii) page 214 of \cite{Kyp}) \\For any $x\in \R$ and $r\geq 0$,
             $$\esp\left(e^{-r\tau^-_0}\1_{\tau^-_0<\infty}|X_0=x\right)=Z^{(r)}(x)-\frac{r}{\Phi(r)}W^{(r)}(x).$$
             \item  It is enough to remark that $\tau^-_x=inf\{t> 0 : X_t<x\}=inf\{t> 0 : 
             X_t-x<0\}$ to deduce that for  $r\geq 0$, the Laplace transform of $\tau^-_x$ is equal to 
             \begin{equation}
             \label{corK4} \esp\left(e^{-r\tau^-_x}\1_{\tau^-_x<\infty}|X_0=0\right)=
             \esp\left(e^{-r\tau^-_0}\1_{\tau^-_0<\infty}|X_0=-x\right)=Z^{(r)}(-x)-\frac{r}{\Phi(r)}W^{(r)}(-x).
             \end{equation}
\end{enumerate}
             \end{proposition}
            
             \begin{preuve} of Proposition \ref{egalite_tau}
             \\
             $\bullet$ {\it{Case x=0}}
             \\
             Remark that $\tau_0=0$ $~\pr(.|X_0=0)$-almost surely.
            We denote   $\tau'_{x}=inf \{t\geq 0 : ~X_t < x\}.$  We want to show that $\tau^-_{0}=0$ and $\tau'_{0}=0$ 
             $~\pr(.|X_0=0)$-almost surely.
             
               Using Proposition \ref{trlap} for $x=0$, we obtain : 
             $$\esp\left(e^{-r\tau^-_0}\1_{\tau^-_0<\infty}|X_0=0\right)=Z^{(r)}(0)-\frac{r}{\Phi(r)}W^{(r)}(0).$$ Since $X$ has unbounded variation, then by Proposition \ref{pptes_w} $(2)$, 
             $W^{(r)}(0)=0$ and by Remark \ref{valen0} $ii)$, $Z^{(r)}(0)=1$, 
             so $$\esp\left(e^{-r\tau^-_0}\1_{\tau^-_0<\infty}|X_0=0\right)=1.$$ However $0\leq e^{-r\tau^-_0}\1_{\tau^-_0<\infty}\leq 1$ implies $e^{-r\tau^-_0}\1_{\tau^-_0<\infty}=1$ $~\pr(.|X_0=0)$-almost surely, and $\tau^-_{0}=0$.
             
              The  process $X$ is smaller than a Brownian motion with drift,  
             $X_t\leq m t+\sigma W_t$ for every $t\geq 0$, and
             $$\pr(\tau'_{0}=0|X_0=0)\geq \pr(\tau^{'m, W}_{x}=0|X_0=0)$$ where  $\tau^{'m, W}_{x}=inf \{t\geq 0 : ~m t+\sigma W_t <  x\}.$ 
             However  $\pr(.|X_0=0)$-almost surely we have  $$inf\{t\geq 0 : m t+\sigma 
             W_t<0\}=inf\{t\geq 0 : m t+\sigma W_t=0\},$$ so $\pr(\tau^{'m, 
             W}_{0}=0|X_0=0)=1$,  consequently $\pr(\tau'_{0}=0|X_0=0)=1.$ 
            
            \vspace{0.3cm}
            
            $ $
            \\
             $\bullet$ {\it{Case x<0}}
             \\
            Since $x<0$, then $\pr(\tau'_{x}=0|X_0=0)=0$ and $\tau'_{x}=\tau^-_{x}$ $~\pr(.|X_0=0)$-almost surely. 
            We want to show that $\tau_{x}=\tau'_{x}$ 
             $~\pr(.|X_0=0)$-almost surely. Since $\tau_x\leq\tau'_{x}$, then 
             $\tau'_{x}=\tau_{x}+\tau'_{x}\circ\theta_{\tau_{x}}$ where $\theta$ is
the translation operator. Apply the strong Markov property at $\tau_{x}$ :
\begin{align*}
\pr(\tau_{x}=\tau'_{x}|X_0=0)&=\esp\left[\esp(\1_{\tau'_{x}=0}|X_0=X_{\tau_{x}})|X_0=0\right]\\
             =&\esp(\1_{X_{\tau_{x}}<x}|X_0=0)+\esp\left(\1_{X_{\tau_{x}}=x}\pr(\tau'_{x}=0|X_0=x)|X_0=0\right).
             \end{align*}
             Since $\pr(\tau'_{x}=0|X_0=x)=\pr(\tau'_{0}=0|X_0=0)=1$, then
             \begin{center}
             $\pr(\tau_{x}=\tau'_{x}|X_0=0)=\pr(X_{\tau_{x}}<x|X_0=0)+\pr(X_{\tau_{x}}=x|X_0=0)=1,$
             \end{center}
             and the conclusion holds.
            \end{preuve}
             \hfill  $\Box$
             \\

             
             According to Proposition \ref{egalite_tau}, the optimal stopping time has the form
 $\tau^-_{ln  \frac{b_c}{v}}$ \\$\pr(.|X_0=0)$-almost surely. We can use  (\ref{corK4}) from Proposition \ref{trlap} and Proposition  \ref{K4} below to find the functions $\cal{L}$ and $\cal{G}$. 
              \begin{proposition} (Exercice 8.7 de \cite{Kyp})
             \label{K4}

              Let $X$ be a spectrally negative process and $\tau^-_x=inf\{t> 
             0 : X_t<x\}$. For any $x<0$, $c\geq 0$ and $r\geq \psi(c)\vee 0$, we have :       $$\esp\left(e^{-r\tau^-_x+c(X_{\tau^-_x}-x)}\1_{\tau^-_x<\infty}|X_0=0\right)=e^{-cx}\left(Z_c^{(r-\psi(c))}(-x)-\frac{r-\psi(c)}{\Phi(r)-c}W_c^{(r-\psi(c))}(-x)\right).$$
             \end{proposition}

\vspace{1cm}

 By Lemma \ref{lem2}, the process $t\mapsto e^{-rt} \left(\frac{-\alpha ve^{X_t}}{r-\psi(1)}+\frac{c}{r}\right)$
 converges in  $L^1$ and almost surely to $0$. Thus, in our case we can remove the indicator function  $\1_{\tau^-_x<\infty}$. The function $g(., b)$ has the form $-\frac{\alpha v}{r-\psi(1)}+ \frac{c}{r}$ if $v\leq b$ ; if not $\frac{-\alpha v}{r-\psi(1)}{\cal{G}}\left(ln\frac{b}{v}\right)+\frac{c}{r}{\cal{L}}\left(ln\frac{b}{v}\right),$
 where 
 \begin{align*}
  {\cal{G}}(x)=&Z_1^{(r-\psi(1))}(-x)-\frac{r-\psi(1)}{\Phi(r)-1}W_1^{(r-\psi(1))}(-x), \\
   {\cal{L}}(x)=&Z^{(r)}(-x)-\frac{r}{\Phi(r)}W^{(r)}(-x).
   \end{align*}
 
 Remark that the function $x\mapsto {\cal{G}}(x)$ is continuous at $x=0$. Indeed, using Definition \ref{def1} and Proposition \ref{pptes_w}, $lim_{x\uparrow 0} {\cal{G}}(x)=1={\cal{G}}(0)$.

In order to be able to check the assumptions of Theorem \ref{th1}, we prove the following results :
            
             \begin{proposition}
             \label{derivees}
             For any $x> 0$, $r\geq \psi(c)\vee 0$ and $c\in \R$, the following equalities are true :
             \begin{enumerate}
             \item
              $Z^{'(r)}(x)=rW^{(r)}(x),$
             \item
             $W_c^{'(r-\psi(c))}(x)=-ce^{-cx}W^{(r)}(x)+e^{-cx}W^{'(r)}(x),$
             \item
             $Z_c^{'(r-\psi(c))}(x)=(r-\psi(c))W_c^{(r-\psi(c))}(x)=(r-\psi(c))e^{-cx}W^{(r)}(x).$
             \end{enumerate}
             \end{proposition}
             \begin{preuve}
             \begin{enumerate}
             \item
              By definition $Z^{(r)}(x)=1+r\int_0^xW^{(r)}(y)dy$, thus
             $Z^{'(r)}(x)=rW^{(r)}(x).$
             \item
             For the second equality, it is enough to differentiate the relation 
             $W_c^{(r-\psi(c))}(x)=e^{-cx}W^{(r)}(x)$. 
             \item
              For the last relation, it is enough to differentiate the function  
             $$x\mapsto Z^{(r-\psi(c))}_c(x)=1+(r-\psi(c))\int_0^xW^{(r-\psi(c))}_c(y)dy$$ and 
             $Z_c^{'(r-\psi(c))}(x)=(r-\psi(c))W_c^{(r-\psi(c))}(x)=(r-\psi(c))e^{-cx}W^{(r)}(x).$
             \end{enumerate}
              \end{preuve}
             \hfill  $\Box$
             \\ 

             The following result is obtained by using Proposition
             \ref{pptes_w} ((2) and (3)) and Proposition \ref{derivees}.
             \begin{corollaire}
             \label{deriv0}
             Let $X$ be a spectrally negative process with unbounded variation. For any $r\geq 
             \psi(c)\vee 0$ and $c\in \R$ :
             \begin{enumerate}
             \item
             $Z^{'(r)}(0)=0$ and $Z_c^{'(r-\psi(c))}(0)=0$,
             \item
              $W_c^{'(r-\psi(c))}(0)=\frac{2}{\sigma^2}$.
              \end{enumerate}
             \end{corollaire}
             \begin{proposition} (proof of  Theorem 9.11 page 258 of \cite{Kyp})
             \label{K5}
             \\
             For any $r\geq 0$ and $x\geq 0$ :
             $$W^{'(r)}(x)-\psi(r)W^{(r)}(x)> 0.$$
             \end{proposition}

              Now we have all the necessary tools to check the assumptions of  Theorem \ref{th1} :
             \\
             \\
               (1) $\cal{G}$ has left derivative at $x=0$.
\\
(2) Moreover \begin{align*}
{\cal{G}}'(0^-)=&-(r-\psi(1))W^{(r)}(0)+\frac{r-\psi(1)}{\Phi(r)-1}\left[ -W^{(r)}(0)+W^{'(r)}(0)\right]\\
=&\frac{r-\psi(1)}{\Phi(r)-1}W^{'(r)}(0)>0
\end{align*}
by Proposition  \ref{K5} for $x=0$. Since ${\cal{L}}'(0^-)=\frac{r}{\Phi(r)}W^{'(r)}(0),$ then 
             $\tilde b=\frac{c(\Phi(r)-1)}{\alpha 
             \Phi(r)}$.  
             \\
             (3) Using Proposition \ref{pptes_w} (4), $g(., \tilde b)\in \cal{C}$$^2(]\tilde b, ~\infty[)$. We use Definition \ref{def1} (2) and Proposition \ref{derivees} to compute the first derivative of $g(., \tilde b)$.

              $\frac{\partial g}{\partial v}(v, \tilde b) =\frac{-\alpha}{r-\psi(1)}Z_1^{(r-\psi(1))}\left(ln \frac{v}{\tilde b}\right)+W^{'(r)}\left(ln \frac{v}{\tilde b}\right)\left( \frac{\alpha \tilde b}{v(\Phi(r)-1)}-\frac{c}{v\Phi(r)}\right)+W^{(r)}\left(ln \frac{v}{\tilde b}\right)\left( \frac{-\alpha \tilde b}{v}+\frac{c}{v}\right).$

             However $\tilde b=\frac{c(\psi(r)-1)}{\alpha \psi(r)}$, so $\frac{\alpha \tilde b}{v(\Phi(r)-1)}-\frac{c}{v\Phi(r)}=0$ and $\frac{-\alpha \tilde b}{v}+\frac{c}{v}=\frac{c}{v\Phi(r)}.$
        $$\frac{\partial g}{\partial 
             v}(v, \tilde b)=-\frac{\alpha}{r-\psi(1)}Z_1^{(r-\psi(1))}\left(ln 
             \frac{v}{\tilde b}\right)+\frac{c}{v\Phi(r)}W^{(r)}\left(ln \frac{v}{\tilde b}\right).$$
             
             By differentiating this relation, we obtain the second derivative of  $g(., \tilde b)$ :
             $$\frac{\partial g^2}{\partial v^2}(v, \tilde b)= -\frac{\alpha}{v(r-\psi(1))} Z_1^{'(r-\psi(1))} \left(ln\frac{v}{\tilde b}\right)-\frac{c}{v^2\Phi(r)}W^{(r)}\left(ln \frac{v}{\tilde b}\right) +\frac{c}{v^2\Phi(r)} W^{'(r)}\left(ln \frac{v}{\tilde b}\right).$$
             
             Using Proposition \ref{derivees},  and replacing $\tilde b$ by its value, we obtain 
            \begin{align*}
            \frac{\partial g^2}{\partial 
             v^2}(v, \tilde b)=&\frac{c}{v^2\Phi(r)}W^{'(r)}\left(ln \frac{v}{\tilde b}\right)-\frac{c}{v^2}W^{(r)}\left(ln \frac{v}{\tilde b}\right)\\
             =&\frac{c}{v^2\Phi(r)}\left[W^{'(r)}\left(ln \frac{v}{\tilde b}\right)-\Phi(r)W^{(r)}\left(ln 
             \frac{v}{\tilde b}\right)\right]>0
             \end{align*}
             by Proposition \ref{K5}. Thus, the function $g(., \tilde b)$ is strictly convex on  $]\tilde b, ~\infty[$.

             We apply Theorem \ref{th1}, $B_c=\tilde b$ and the smallest optimal stopping time is 
             $$\tau^*(c)=inf \{t\geq 0 : ~V_t \leq B_c\}.$$

             \begin{proposition}
             Let $X$ be a spectrally negative process with a non null Gaussian component.  Then, with the notations introduced in Section \ref{levneg},
             \begin{enumerate}
             \item
           The smallest optimal stopping time is
             $\tau^*(c)=inf \{t\geq 0 : ~V_t \leq B_c\}$ where \\$B_c=\frac{c(\Phi(r)-1)}{\alpha 
             \Phi(r)}$.
             \item
              For $v>~B_c$, the value function  $w$ is equal to \\$\frac{\alpha v}{r-\psi(1)}-\frac{c}{r}+\frac{-\alpha 
             v}{r-\psi(1)}\left[Z_1^{(r-\psi(1))}\left(ln 
             \frac{v}{B_c}\right)-\frac{r-\psi(1)}{\psi(r)-1}W_1^{(r-\psi(1))}\left(ln 
             \frac{v}{B_c}\right)\right]+\frac{c}{r}\left[Z^{(r)}\left(ln 
             \frac{v}{B_c}\right)-\frac{r}{\psi(r)}W^{(r)}\left(ln \frac{v}{B_c}\right)\right].$
             \end{enumerate}
            \end{proposition}

            \vspace{1cm}
            
            \textbf{Conclusion}
            
            Our method is much easier than the traditional methods (Wiener-Hopf factorization -- see for exemple  \cite{boy1, boy2, Kou2},  Monte-Carlo method as in \cite{Long} or  integro-differential equations as in \cite{Duf, Mord}). It can be used for all Lévy process when the the joint  Laplace transform of  $(\tau_b, ~X_{\tau_b})$ is known, where $\tau_b=inf\{t\geq 0 : ~X_t \leq b\}$. 
             
             \section{Appendix}
             \label{appendices}
             \subsection{ Optimal stopping tools}
             \label{appendix1}
             
             For the sake of completeness, we recall some classical results of optimal stopping theory  used to solve the problem studied in this paper (we refer to \cite{KLM} and \cite{Sh}).

             \begin{theoreme} (Theorem 3.4  of  \cite{KLM})
             \label{formeSnell}

             Let $V_.$ be a strong Markov process and $Y_.$ a process of class $D$
            of the form $t\mapsto e^{-rt}f(V_t)$ where $f$ is a  measurable function.
             Let $J$ be its  Snell envelope (i.e. the smallest supermartingale larger than $Y$) :  $J_t=esssup_{\tau\in\Delta , \tau \geq t}E[Y_\tau \mid\F_t]$.
              Then $J_.$ has the form $t\mapsto J_t=e^{-rt}s(V_t)$ where the function $s$ is called "r-reduite" of $f$.
             \end{theoreme}

             \begin{theoreme} (Optimality criteria -- Remark 3.5 of \cite{KLM})
             \label{toptim}

               Let $Y_.$ be a strong Markov process  and  $J$ its  Snell envelope. A stopping time $\tau^*$ is optimal if and only if :

             -  $Y_{\tau^*}= J_{\tau^*}$,

             -  $J_{. \wedge\tau^*}$ is a martingale.

             \end{theoreme}

             \begin {theoreme}  (Theorem 3.3 page 127  of \cite{Sh})
             \label{th5}
             
             Let $V_.$ be a strong Markov process and $Y_.$ a process of the form $t\mapsto f(V_t)$ where $f$ is a measurable function.  Let $J$ be  the Snell envelope of $Y$. 
             
             For any $\eps\geq 0$, let $\tau_{\eps}=inf\{t\geq 0 : J_t\leq Y_t +\eps\}$. If $Y_t$ satisfies the following conditions : 
             $\pr(lim_{t \downarrow 0} Y_t=Y_0)=1$, ~$\esp[sup_{t\geq 0} max(Y_t, 
             0)]<\infty$,
             ~$\esp[sup_{t\geq 0} -min(Y_t, 0)]<\infty$, then

             1. For any $\eps>0$, the times $\tau_{\eps}$ are $\eps$-optimal stopping times.

             2. If the function $f$ is upper semicontinuous, i.e. $\overline{lim}_{y\rightarrow x} 
             f(y) \leq f(x)$, then $\tau_0$ is an optimal stopping time.

             3. If there exists an optimal time $\tau \in  \Delta$, then
             $\pr(\tau_0\leq \tau)=1$ and $\tau_0 \in  \Delta$  is optimal.
             \end{theoreme}
             
             \begin{remarque}
             i) Under the hypothesis of Theorem \ref{th5}, if the function $f$ is u.s.c., then using 3., $\tau_0$ is the smallest optimal stopping time.
             \\
             ii) Theorem \ref{th5} is also checked when $Y$ is a process of the form  $t\mapsto e^{-rt}\bar{f}(V_t)$ where $\bar{f}$ is a measurable function. Indeed, since  $X_t=(V_t, ~t), ~t\geq 0$ is a strong Markov process it is enough to consider $f(X_t)=e^{-rt}\bar{f}(V_t)$.
             \end{remarque}

             \begin{lemme} (Lemma 3.8. page 123 of \cite{Sh})
             \label{tpepsopt}

             Let $Y$ be a strong Markov process and $J$ its Snell envelope. For any $\eps\geq 0$, let \\$\tau_{\eps}=inf\{t\geq 0 : J_t\leq Y_t +\eps\}$.
             If the process $Y$ satisfies the following conditions \\$\pr(lim_{t \downarrow 0} Y_t=Y_0)=1$ and $\esp[sup_{t\geq 0} max(Y_t, 0)]<\infty$, then for any $\eps>0$, 
             $P(\tau_{\eps}<\infty)=1$.
             \end{lemme}
             
  
             
             \begin{theoreme} (Theorem 25 page 92 of \cite{DM})
             \label{classeD_DM}
             
            A positive right continuous supermartingale $X$ is of class $D$ if and only if  $$lim_{n\rightarrow\infty}\esp\left(X_{R_n}\1_{R_n<\infty}\right)=0$$ where $R_n=inf\{t\geq 0 : X_t\geq n\}.$
            \end{theoreme}

             \subsection{Useful result}
             \label{appendix2}
             
             Next, we present a useful result for the calculation of the optimal strategy in the case of a  particular mixed  diffusion-jump process. The following proposition starts from \cite{Kou} and \cite{Dao}.
In  \cite{Kou} and \cite{Dao}, the authors calculate the Laplace transform of a first passage time of the form
 $\tau^b=inf\{t\geq 0 :  X_t \geq b\}$ where $X$ is a mixed diffusion-jump process and the jump size is a random variable with a double exponential distribution. They also calculate the joint Laplace transform of  $(\tau^b, ~X_{\tau^b})$  and give the calculation algorithm in the case of a first passage time of the form $\tau_b=inf\{t\geq 0 : X_t \leq b\}$.

Throughout this section   $(W_t, ~t\geq 0)$ is a standard Brownian motion, $(N_t, ~t\geq 0)$  a Poisson process with constant positive intensity $a$, $(Y_i, ~i\in \N )$ is a sequence of independent and identically distributed random variables with a double exponential distribution, i.e. the common density of $Y$ is given by
             $$f_Y(y)=p\eta_1e^{-\eta_1y}1_{y>0}+q\eta_2e^{\eta_2y}1_{y<0}$$
             where $p+q=1$, $p,q>0$, $\eta_1>1$ and $\eta_2>0$.
             \begin{proposition} 
             \label{kou}

             Let $\tau$ be a the first passage time of the form  $$\tau=inf\{t\geq 0 : 
           X_t \leq b\}$$ where $m\in\R$, $\sigma>0$, $b<0$ and $X_t= mt+\sigma W_t+\sum_{i=1}^{N_t} Y_i, ~t\geq 0$. For any $r\geq 0$ :

             \begin{equation}
             \label{14}
             \esp[e^{-r\tau}|X_0=0]=\frac{\psi_2(\eta_2+\psi_3)}{(\psi_2-\psi_3)\eta_2}e^{-b\psi_3}-\frac{\psi_3(\eta_2+\psi_2)}{(\psi_2-\psi_3)\eta_2}e^{-b\psi_2},
             \end{equation}

             \begin{equation}
             \label{15}
             \esp[e^{-r\tau+X_{\tau}}\1_{\tau<\infty}|X_0=0]=e^{b}\left[ \frac 
             {(\eta_2+\psi_3)(\psi_2-1)}{(\psi_2-\psi_3)(\eta_2+1)}e^{-b\psi_3}+ 
             \frac{(\eta_2+\psi_2)(1-\psi_3)}{(\psi_2-\psi_3)(\eta_2+1)}e^{-b\psi_2}\right],
             \end{equation}
            where $\psi_2$, $\psi_3$ are the negative roots of the equation 
            \begin{equation}
            \label{exponentlevy}
             m\psi+\frac{\sigma^2}{2}\psi^2+a[\frac{\eta_1p}{\eta_1-\psi}+\frac{\eta_2q}{\eta_2+\psi}-1]=r,
             \end{equation}
              $-\infty<\psi_3<-\eta_2<\psi_2<0$. 
             \end{proposition}

             \begin{preuve}
\\
             Using  Theorem 3.1 of \cite{Kou}, the Laplace transform of the following stopping time \\ $\tau^b=inf\{t\geq 0 : mt+\sigma  W_t+\sum_{i=1}^{N_t} Y_i \geq b\}$ where $b>0$, is :    $$\esp[e^{-r\tau^b}|X_0=0]=\frac{\psi_0(\eta_1-\psi_1)}{(\psi_0-\psi_1)\eta_1}e^{-b\psi_1}+\frac{\psi_1(-\eta_1+\psi_0)}{(\psi_0-\psi_1)\eta_1}e^{-b\psi_0}$$

             where $0<\psi_1<\eta_1<\psi_0<\infty$ are the positive roots of the equation
             \begin{equation}
             \nonumber
             m\psi+\frac{\sigma^2}{2}\psi^2+a[\frac{\eta_1p}{\eta_1-\psi}+\frac{\eta_2q}{\eta_2+\psi}-1]=r.
             \end{equation}

             By Corollary 3.3  of \cite{Kou}, for any $\beta<\eta_1$, 
\begin{equation}
\label{1_ch4}
             \esp[e^{-r\tau^b+\beta X_{\tau^b}}\1_{\tau^b<\infty}|X_0=0]=e^{\beta b}[ \frac 
             {(\eta_1-\psi_1)(\psi_0-\beta)}{(\psi_0-\psi_1)(\eta_1-\beta)}e^{-b\psi_1}+ 
             \frac{(\psi_0-\eta_1)(\psi_1-\beta)}{(\psi_0-\psi_1)(\eta_1-\beta)}e^{-b\psi_0}].
             \end{equation}

             The result follows from Remark 4.2 \cite{Kou}. Indeed, according to Remark 4.2 of \cite{Kou} or \cite{Dao}, to make the same calculation as above
            for a stopping time of the form \\
             $\tau=inf\{t\geq 0 : mt+\sigma W_t+\sum_{i=1}^{N_t} Y_i \leq b\}$ where 
             $b<0$, we only need to make the following changes : $p\mapsto q$, 
             $q\mapsto p$, $\psi_1\mapsto -\psi_2$, $\psi_0\mapsto -\psi_3$, 
             $\eta_1\mapsto \eta_2$, $\eta_2\mapsto \eta_1$, $b\mapsto -b$, $\beta\mapsto -\beta$ 
             where $\psi_2$ and $\psi_3$ are the negative roots of the equation  (\ref{exponentlevy}). 
              \end{preuve}
             \hfill  $\Box$
             \\

\begin{remarque}
\label{p1}
Even if $p=1$ (and thus $q=0$), the relation (\ref{1_ch4}) is true. But, since $p=1$, the equation (\ref{exponentlevy}) has only one negative root and in this case (\ref{14}) and (\ref{15}) are not true.
\end{remarque}

             \end{document}